\documentclass[12pt]{article}
\usepackage{amsfonts}
\usepackage{full page, 
amssymb, amscd, amsmath,graphicx, 
}

 \title{{\bf Intertwining operators among twisted modules associated to 
not-necessarily-commuting automorphisms}}
 \author{Yi-Zhi Huang}
   \date{}
    \begin{document}
    \bibliographystyle{alpha}
    \maketitle
\newtheorem{thm}{Theorem}[section]
\newtheorem{defn}[thm]{Definition}
\newtheorem{prop}[thm]{Proposition}
\newtheorem{cor}[thm]{Corollary}
\newtheorem{lemma}[thm]{Lemma}
\newtheorem{rema}[thm]{Remark}
\newtheorem{app}[thm]{Application}
\newtheorem{prob}[thm]{Problem}
\newtheorem{conv}[thm]{Convention}
\newtheorem{conj}[thm]{Conjecture}
\newtheorem{cond}[thm]{Condition}
    \newtheorem{exam}[thm]{Example}
\newtheorem{assum}[thm]{Assumption}
     \newtheorem{nota}[thm]{Notation}
\newcommand{\halmos}{\rule{1ex}{1.4ex}}
\newcommand{\pfbox}{\hspace*{\fill}\mbox{$\halmos$}}

\newcommand{\nn}{\nonumber \\}

 \newcommand{\res}{\mbox{\rm Res}}
 \newcommand{\ord}{\mbox{\rm ord}}
\renewcommand{\hom}{\mbox{\rm Hom}}
\newcommand{\edo}{\mbox{\rm End}\ }
 \newcommand{\pf}{{\it Proof.}\hspace{2ex}}
 \newcommand{\epf}{\hspace*{\fill}\mbox{$\halmos$}}
 \newcommand{\epfv}{\hspace*{\fill}\mbox{$\halmos$}\vspace{1em}}
 \newcommand{\epfe}{\hspace{2em}\halmos}
\newcommand{\nord}{\mbox{\scriptsize ${\circ\atop\circ}$}}
\newcommand{\wt}{\mbox{\rm wt}\ }
\newcommand{\swt}{\mbox{\rm {\scriptsize wt}}\ }
\newcommand{\lwt}{\mbox{\rm wt}^{L}\;}
\newcommand{\rwt}{\mbox{\rm wt}^{R}\;}
\newcommand{\slwt}{\mbox{\rm {\scriptsize wt}}^{L}\,}
\newcommand{\srwt}{\mbox{\rm {\scriptsize wt}}^{R}\,}
\newcommand{\clr}{\mbox{\rm clr}\ }
\newcommand{\tr}{\mbox{\rm Tr}}
\newcommand{\C}{\mathbb{C}}
\newcommand{\Z}{\mathbb{Z}}
\newcommand{\R}{\mathbb{R}}
\newcommand{\Q}{\mathbb{Q}}
\newcommand{\N}{\mathbb{N}}
\newcommand{\CN}{\mathcal{N}}
\newcommand{\F}{\mathcal{F}}
\newcommand{\I}{\mathcal{I}}
\newcommand{\V}{\mathcal{V}}
\newcommand{\one}{\mathbf{1}}
\renewcommand{\i}{\mathbf{i}}
\newcommand{\BY}{\mathbb{Y}}
\newcommand{\ds}{\displaystyle}

        \newcommand{\ba}{\begin{array}}
        \newcommand{\ea}{\end{array}}
        \newcommand{\be}{\begin{equation}}
        \newcommand{\ee}{\end{equation}}
        \newcommand{\bea}{\begin{eqnarray}}
        \newcommand{\eea}{\end{eqnarray}}
         \newcommand{\lbar}{\bigg\vert}
        \newcommand{\p}{\partial}
        \newcommand{\dps}{\displaystyle}
        \newcommand{\bra}{\langle}
        \newcommand{\ket}{\rangle}

        \newcommand{\ob}{{\rm ob}\,}
        \renewcommand{\hom}{{\rm Hom}}

\newcommand{\A}{\mathcal{A}}
\newcommand{\Y}{\mathcal{Y}}

\newcommand{\dlt}[3]{#1 ^{-1}\delta \bigg( \frac{#2 #3 }{#1 }\bigg) }

\newcommand{\dlti}[3]{#1 \delta \bigg( \frac{#2 #3 }{#1 ^{-1}}\bigg) }

 \makeatletter
\newlength{\@pxlwd} \newlength{\@rulewd} \newlength{\@pxlht}
\catcode`.=\active \catcode`B=\active \catcode`:=\active
\catcode`|=\active
\def\sprite#1(#2,#3)[#4,#5]{
   \edef\@sprbox{\expandafter\@cdr\string#1\@nil @box}
   \expandafter\newsavebox\csname\@sprbox\endcsname
   \edef#1{\expandafter\usebox\csname\@sprbox\endcsname}
   \expandafter\setbox\csname\@sprbox\endcsname =\hbox\bgroup
   \vbox\bgroup
  \catcode`.=\active\catcode`B=\active\catcode`:=\active\catcode`|=\active
      \@pxlwd=#4 \divide\@pxlwd by #3 \@rulewd=\@pxlwd
      \@pxlht=#5 \divide\@pxlht by #2
      \def .{\hskip \@pxlwd \ignorespaces}
      \def B{\@ifnextchar B{\advance\@rulewd by \@pxlwd}{\vrule
         height \@pxlht width \@rulewd depth 0 pt \@rulewd=\@pxlwd}}
      \def :{\hbox\bgroup\vrule height \@pxlht width 0pt depth
0pt\ignorespaces}
      \def |{\vrule height \@pxlht width 0pt depth 0pt\egroup
         \prevdepth= -1000 pt}
   }
\def\endsprite{\egroup\egroup}
\catcode`.=12 \catcode`B=11 \catcode`:=12 \catcode`|=12\relax
\makeatother

\def\hboxtr{\FormOfHboxtr} 
\sprite{\FormOfHboxtr}(25,25)[0.5 em, 1.2 ex] 

:BBBBBBBBBBBBBBBBBBBBBBBBB | :BB......................B |
:B.B.....................B | :B..B....................B |
:B...B...................B | :B....B..................B |
:B.....B.................B | :B......B................B |
:B.......B...............B | :B........B..............B |
:B.........B.............B | :B..........B............B |
:B...........B...........B | :B............B..........B |
:B.............B.........B | :B..............B........B |
:B...............B.......B | :B................B......B |
:B.................B.....B | :B..................B....B |
:B...................B...B | :B....................B..B |
:B.....................B.B | :B......................BB |
:BBBBBBBBBBBBBBBBBBBBBBBBB |

\endsprite
\def\shboxtr{\FormOfShboxtr} 
\sprite{\FormOfShboxtr}(25,25)[0.3 em, 0.72 ex] 

:BBBBBBBBBBBBBBBBBBBBBBBBB | :BB......................B |
:B.B.....................B | :B..B....................B |
:B...B...................B | :B....B..................B |
:B.....B.................B | :B......B................B |
:B.......B...............B | :B........B..............B |
:B.........B.............B | :B..........B............B |
:B...........B...........B | :B............B..........B |
:B.............B.........B | :B..............B........B |
:B...............B.......B | :B................B......B |
:B.................B.....B | :B..................B....B |
:B...................B...B | :B....................B..B |
:B.....................B.B | :B......................BB |
:BBBBBBBBBBBBBBBBBBBBBBBBB |

\endsprite

\vspace{2em}



\renewcommand{\theequation}{\thesection.\arabic{equation}}
\renewcommand{\thethm}{\thesection.\arabic{thm}}
\date{}
\maketitle

\begin{abstract}
We introduce intertwining operators among twisted modules or twisted intertwining operators
associated to not-necessarily-commuting automorphisms
of a vertex operator algebra. Let $V$ be a vertex operator algebra and let $g_{1}$, $g_{2}$ and $g_{3}$ be automorphisms 
of $V$. 
We prove that for $g_{1}$-, $g_{2}$- and  $g_{3}$-twisted $V$-modules 
$W_{1}$, $W_{2}$ and $W_{3}$, respectively, such that the vertex operator map for $W_{3}$ is injective,
if there exists a twisted intertwining 
operator of type ${W_{3}\choose W_{1}W_{2}}$ such that the images of its component
operators span $W_{3}$, then $g_{3}=g_{1}g_{2}$.
We also construct what we call the skew-symmetry 
and contragredient isomorphisms between spaces of twisted intertwining operators among twisted modules of suitable types. 
The proofs of these results involve careful
analysis of the analytic extensions corresponding to 
the actions of the not-necessarily-commuting automorphisms of the vertex operator algebra. 
\end{abstract}

\setcounter{equation}{0} \setcounter{thm}{0}

\section{Introduction}

In the present paper, we initiate the study of intertwining operators among twisted modules 
associated to not-necessarily-commuting automorphisms of a vertex operator algebra. 
Here by twisted modules we mean (generalized or logarithmic) twisted modules introduced in \cite{H8}. For simplicity
and to avoid confusion, when twisted modules are not mentioned, 
we shall call such intertwining operators 
``twisted intertwining operators,'' although these intertwining 
operators are not twisted directly, but are twisted instead in a suitable sense through twisted modules.  

Intertwining operators
among (untwisted) modules for a vertex operator algebra
were first introduced mathematically by Frenkel, Lepowsky and the author
in \cite{FHL} and  correspond to chiral vertex operators in physics 
(see \cite{MS}).  They have been studied systematically in the papers \cite{HL}, \cite{H1}--\cite{H6}, 
\cite{HLZ2}--\cite{HLZ7},  \cite{Y}, \cite{C1}--\cite{C2} and  \cite{F1}--\cite{F2}.  Intertwining operators give
chiral three-point correlation functions in two-dimensional conformal field theories and are
the building blocks of multi-point correlation functions on Riemann surfaces of arbitrary genus. They
are the main objects of interest in the representation theory of vertex operator algebras and 
two-dimensional conformal field theory. Almost all important results in these theories are in fact
properties of  intertwining operators.

Intertwining operators among twisted modules associated 
to commuting automorphisms of finite orders appeared implicitly in the work \cite{FFR}  of  Feingold, Frenkel and Ries 
and were introduced explicitly by Xu in \cite{X} in terms of  a generalization of the Jacobi identity for twisted modules. 
However, there is still no definition of intertwining operators among twisted modules associated to noncommuting 
automorphisms in the literature. To construct orbifold conformal field theories associated to a noncommutative
group of automorphisms of a vertex operator algebra, it is necessary to study these intertwining operators. 

In \cite{H9}, the author formulated the following conjecture:

\begin{conj}[\cite{H9}]\label{twisted-int-op}
Assume that $V$ is a simple vertex operator algebra 
satisfying the following conditions:
\begin{enumerate}
\item $V_{(0)}=\C\one$, $V_{(n)}=0$ 
for $n<0$ and the contragredient $V'$, as a $V$-module, is equivalent to $V$.

\item Every grading-restricted generalized $V$-module is completely reducible

\item $V$ is $C_{2}$-cofinite, that is, $\dim V/C_{2}(V)<\infty$, where $C_{2}(V)$ is the subspace of 
$V$ spanned by the elements of the form $\res_{x}x^{-2}Y(u, x)v$ for $u, v\in V$ and $Y: V\otimes V\to V[[x, x^{-1}]]$
is the vertex operator map for $V$.
\end{enumerate}
Let $G$ be a finite group 
of automorphisms of $V$. Then the twisted intertwining operators among the $g$-twisted $V$-modules for all $g\in G$
satisfy the associativity, commutativity and modular invariance properties.
\end{conj}

\noindent If this conjecture is proved, then we obtain the genus-zero and genus-one parts of the chiral
orbifold conformal field theory associated with the vertex operator algebra $V$ and the group
$G$ of automorphisms of $V$. 
One consequence of Conjecture \ref{twisted-int-op} is that the category of $g$-twisted modules for all 
all $g\in G$ has a natural structure of $G$-crossed braided tensor category satisfying additional properties. 

To even formulate this conjecture precisely, we have to first introduce twisted intertwining operators or intertwining operators
among twisted modules associated to general automorphisms and study their basic properties. 
In this paper, we give a definition of such twisted intertwining operators. 
Let $V$ be a vertex operator algebra and let $g_{1}$, $g_{2}$ and $g_{3}$ be automorphisms 
of $V$. 
We prove that for $g_{1}$-, $g_{2}$- and  $g_{3}$-twisted $V$-modules 
$W_{1}$, $W_{2}$ and $W_{3}$, respectively, such that the vertex operator map for $W_{3}$ is injective,
if there exists a twisted intertwining 
operator of type ${W_{3}\choose W_{1}W_{2}}$ such that the images of its component
operators span $W_{3}$, then $g_{3}=g_{1}g_{2}$.
We also construct what we call the skew-symmetry 
and contragredient isomorphisms between spaces of twisted intertwining operators among twisted modules of suitable types. 
The proofs of these results are much more subtle and delicate than those for the corresponding results in \cite{FHL}, \cite{HL}, 
\cite{X} and \cite{HLZ2}  because they involve careful
analysis of the analytic extensions corresponding to 
the actions of the not-necessarily-commuting automorphisms of the vertex operator algebra.

The motivation for the study of twisted intertwining operators does not just come from orbifold
conformal field theories and their potential applications in geometry and physics. 
It in fact also comes intrinsically from the study of the uniqueness conjecture
of the moonshine module vertex operator algebra proposed by Frenkel, Lepowsky and Meurman \cite{FLM3}. 
This conjecture was obtained by using the analogy among the Golay code, the Leech lattice and 
the moonshine module vertex operator algebra. In Conway's proof of the uniqueness of the Leech lattice
\cite{Co}, the $24$-dimensional vector space $\R^{24}$  plays a fundamental role. 
For the uniqueness conjecture for the moonshine module vertex operator algebra, the first difficulty 
is that there is no analogue of $\R^{24}$. But we still 
need a structure large enough  such that all the works can be done in this structure. 

If the vertex operator subalgebra  fixed by the automorphism group of a vertex operator algebra
satisfy suitable conditions (for example, the three conditions in Conjecture \ref{twisted-int-op}),  
then one possible choice of such a structure large enough for our purpose 
is the intertwining operator algebra formed by the modules and intertwining operators for the fixed-point vertex operator 
subalgebra. However, the assumption that these suitable conditions hold
is in fact a main difficult conjecture that we have to prove first. Because of this, 
instead of assuming these conditions, we have to develop a 
theory that will lead us to a proof of these conditions and a construction of 
the intertwining operator algebra. 

Since in this area, conjectures are sometimes claimed to have been proved in some books, papers, preprints
or unpublished manuscripts without the evidence that the proofs indeed exist, the author would like to comment that 
the mathematics for obtaining conjectures and the mathematics for finding proofs are often very different. 
To obtain conjectures, one can assume what one believes to be true and derive the consequences. But to
prove conjectures, the assumptions  used to derive the conjectures are often themselves the most difficult parts of the
conjectures. Thus one might need different mathematical approaches and theories to prove these assumptions first.
The theory of twisted intertwining operators initiated in this paper is what the author believes to be needed 
in the proofs of the conjectures mentioned above, including 
Conjecture \ref{twisted-int-op} and the conjecture that suitable conditions hold for the  fixed-point vertex operator 
subalgebra. We expect that this theory will 
play an important role in the study of orbifold conformal field theory and, in particular, 
in the study of the uniqueness conjecture of the moonshine module vertex operator algebra.

In this paper, we use the approach based on multivalued analytic functions with preferred branches developed 
and used in  \cite{H1}, \cite{H3}, \cite{H6}, \cite{HLZ5}--\cite{HLZ6} and \cite{C1}--\cite{C2}. Intertwining operators among (untwisted) modules
can be defined using either such multivalued analytic functions or formal variables. But to study products and iterates
of intertwining operators, it is necessary to use the approach based on such multivalued analytic functions. 
For twisted intertwining operators introduced and studied in this paper, even for the definition and the properties involving only 
one twisted intertwining operator, we need the approach based on such multivalued analytic functions, because the vertex operators
for twisted modules contain nonintegral powers and logarithm of the variable and, more importantly,
because the automorphisms associated to different twisted modules do not necessarily commute with each other. 

In Section 2, we discuss the notations and conventions used in this paper, 
especially those involving multivalued analytic functions with preferred branches. 
In Section 3, we recall the notion of twisted module introduced in \cite{H8}. We also discuss 
in this section the functors
associated to automorphisms of the vertex operator algebra and the contragredient 
functor on the category of twisted modules.  Twisted intertwining 
operators are introduced in Section 4. In the same section, 
we prove the result that for $g_{1}$-, $g_{2}$- and  $g_{3}$-twisted $V$-modules 
$W_{1}$, $W_{2}$ and $W_{3}$, respectively, such that the vertex operator map for $W_{3}$ is injective,
if there exists a  twisted intertwining 
operator of type ${W_{3}\choose W_{1}W_{2}}$ such that the images of its component
operators span $W_{3}$, then $g_{3}=g_{1}g_{2}$.
The skew-symmetry  isomorphisms 
and contragredient isomorphisms are constructed in Section 5 and Section 6, respectively.

\setcounter{equation}{0} \setcounter{thm}{0}

\section{Notations and conventions}

To study intertwining operators, we have to work with multivalued analytic functions with preferred branches. 
The approach that we use in this paper is the same as the one used in \cite{H1}, \cite{H3}, \cite{H6},
\cite{HLZ5}--\cite{HLZ6} and \cite{C1}--\cite{C2}. In this section, we recall and introduce
some notations and conventions. 

We shall use $\i$ to denote $\sqrt{-1}$. 
For $z\in \C^{\times}$, we choose the value $\arg z$ of the argument of $z$ to be 
the one satisfying $0\le \arg z<2\pi$. We shall not use $\log z$ to denote the particular 
value $\log |z|+ (\arg z)\i$ of the logarithm of $z$ as in  \cite{H1}, \cite{H3}, 
\cite{HLZ5}--\cite{HLZ6} and \cite{C1}--\cite{C2}. Instead, we shall always use 
$l_{p}(z)$ to denote the value $\log |z|+(\arg z+ 2p\pi) \i$ of the logarithm of $z$ for
$p\in \Z$. 

Intertwining operators defined using formal variables in fact give multivalued analytic
functions with preferred branches. 
We shall use $\log z$ to denote the multivalued logarithm function of the 
variable $z$ with the preferred branch $l_{0}(z)=\log |z|+(\arg z) \i$. For $n\in \C$, we shall 
use $z^{n}$ to denote the multivalued analytic function $e^{n\log z}$ with the 
preferred branch $e^{nl_{0}(z)}$. Multivalued analytic functions with preferred 
branch on a region form a commutative associative algebra and can also be 
divided by such functions on the same region to obtain such functions on possibly smaller regions. In particular, 
$$f(z_{1}, z_{2})=\sum_{i, j, k, l, m, n=1}^{N}a_{ijklmn}z_1^{r_{i}}z_2^{s_{j}}(z_1 - z_2)^{t_{k}}(\log z_1)^l
(\log z_2)^{m}(\log (z_1 - z_2))^{n}$$
for $a_{ijklmn}, r_{i}, s_{j}, t_{k}\in \C$ is a multivalued analytic function with preferred branch 
on the region given by $z_{1}, z_{2}\ne 0$, $z_{1}\ne z_{2}$. 
For $p_{1}, p_{2}, p_{12}\in \Z$, we shall use 
$f^{p_{1}, p_{2}, p_{12}}(z_{1}, z_{2})$
to denote the single-valued branch 
$$\sum_{i, j, k, l, m, n=1}^{N}a_{ijklmn}e^{r_{i}l_{p_{1}}(z_{1})}
e^{s_{j}l_{p_{2}}(z_2)}e^{t_{k}l_{p_{12}}(z_1 - z_2)}(l_{p_{1}}(z_1))^l
(l_{p_{2}}(z_2))^{m}(l_{p_{12}}(z_1 - z_2))^{n}$$
of $f(z_{1}, z_{2})$. 

For a $\C$-graded vector space $W=\coprod_{n\in \C}W_{[n]}$, let $W'=\coprod_{n\in \C}W_{[n]}^{*}$ be the
graded dual of $W$ and $\overline{W}=\prod_{n\in \C}W_{[n]}$  the algebraic completion of $W$.
For $n \in \mathbb{C}$, we use $\pi_n$ to denote the  the projection from $W$ or $\overline{W}$ to $W_{[n]}$ . 

Let $W$ be a vector space and 
$$X(x)=\sum_{k=0}^{K}\sum_{n\in \C}a_{n, k}x^{n}(\log x)^{k}\in (\text{\rm End}\;W)\{x\}[\log x].$$
For $z\in \C^{\times}$, we shall use 
$X^{p}(z)$ to denote the series
$$X(x)\lbar_{x^{n}=e^{nl_{p}(z)},\; \log x=l_{p}(z)}=\sum_{k=0}^{K}\sum_{n\in \C}a_{n, k}e^{nl_{p}(z)}(l_{p}(z))^{k}$$
with terms in $\text{\rm End}\;W$. When $W=\coprod_{n\in \C}W_{[n]}$ is a $\C$-graded vector space and $a_{n, k}$ 
for different $n$ are homogeneous operators of different degrees, 
$X^{p}(z)\in \hom(W, \overline{W})$. When $z$ changes in $\C^{\times}$, $X^{p}(z)$ can be viewed as a 
function on $\C^{\times}$ valued in the space of series in $W$. 
We call  this function $X^{p}(z)$ the {\it $p$-th analytic branch of $X(x)$}.

\setcounter{equation}{0} \setcounter{thm}{0}

\section{Twisted modules}

In this paper, we fix a vertex operator algebra $(V, Y_{V}, \one_{V}, \omega_{V})$. 
In fact, the results of the present paper hold for 
a grading-restricted M\"{o}bius vertex algebra, that is, a $\Z$-graded vertex algebra $V=\coprod_{n\in \Z}V_{(n)}$
equipped with operators $L_{V}(-1)$, $L_{V}(0)$ and $L_{V}(1)$
such that $V_{(n)}=0$ when $n$ is sufficiently negative, $\dim V_{(n)}<\infty$ for $n\in \Z$,
the usual $\mathfrak{sl}(2, \C)$ commutator relations for $L_{V}(-1)$, $L_{V}(0)$ and $L_{V}(1)$  hold and 
the usual commutator relations between $L_{V}(-1)$, $L_{V}(0)$ and $L_{V}(1)$  and vertex operators
hold. But for what we want to prove in the future, it is necessary for $V$ to be a vertex operator
algebra with the additional data of a conformal element satisfying some additional conditions. 

Let $g$ be an automorphism of $V$.
We first recall the definition of generalized $g$-twisted $V$-module first introduced in \cite{H8}. But for simplicity,
we shall omit the word ``generalized.'' In particular,  in this paper, the vertex operator map for a $g$-twisted $V$-module
in general contain the logarithm of the variable and the operator $L(0)$ in general does not have to act semisimply.

\begin{defn}
{\rm A {\it  $g$-twisted
$V$-module} is a ${\C}\times \C/\Z$-graded
vector space $W = \coprod_{n \in \C, \alpha\in \C/\Z} W_{[n]}^{[\alpha]}$ (graded by {\it weights} and {\it $g$-weights})
equipped with a linear map
\begin{eqnarray*}
Y_{W}^g: V\otimes W &\to& W\{x\}[\log  x],\\
v \otimes w &\mapsto &Y_{W}^g(v, x)w
\end{eqnarray*}
satisfying the following conditions:
\begin{enumerate}

\item The {\it equivariance property}: For $p \in \mathbb{Z}$, $z
\in \mathbb{C}^{\times}$,  $v \in V$ and $w \in W$, 
$$(Y^{g}_{W})^{p + 1}(gv,
z)w = (Y^{g}_{W})^{p}(v, z)w,$$
where for $p \in \mathbb{Z}$, $(Y^{g}_{W})^{ p}(v, z)$
is the $p$-th analytic branch of $Y_{W}^g(v, x)$.

\item The {\it identity property}: For $w \in W$, $Y^g({\bf 1}, x)w
= w$.

\item  The {\it duality property}: For
any $u, v \in V$, $w \in W$ and $w' \in W'$, there exists a
multivalued analytic function with preferred branch of the form
\[
f(z_1, z_2) = \sum_{i,
j, k, l = 0}^N a_{ijkl}z_1^{m_i}z_2^{n_j}({\rm log}z_1)^k({\rm
log}z_2)^l(z_1 - z_2)^{-t}
\]
for $N \in \mathbb{N}$, $m_1, \dots,
m_N$, $n_1, \dots, n_N \in \mathbb{C}$ and $t \in \mathbb{Z}_{+}$,
such that the series
\[
\langle w', (Y^{g}_{W})^{ p}(u, z_1)(Y^{g}_{W})^{p}(v,
z_2)w\rangle = \sum_{n \in \mathbb{C}}\langle w', (Y^{g}_{W})^{ p}(u,
z_1)\pi_n(Y^{g}_{W})^{ p}(v, z_2)w\rangle,
\]
\[
\langle w', (Y^{g}_{W})^{ p}(v,
z_2)(Y^{g}_{W})^{p}(u, z_1)w\rangle = \sum_{n \in \mathbb{C}}\langle w',
(Y^{g}_{W})^{ p}(v, z_2)\pi_n(Y^{g}_{W})^{ p}(u, z_1)w\rangle,
\]
\[
\langle w', (Y^{g}_{W})^{ p}(Y_{V}(u, z_1 - z_2)v,
z_2)w\rangle = \sum_{n \in \mathbb{C}}\langle w', (Y^{g}_{W})^{
p}(\pi_nY_{V}(u, z_1 - z_2)v, z_2)w\rangle
\]
are absolutely convergent in
the regions $|z_1| > |z_2| > 0$, $|z_2| > |z_1| > 0$, $|z_2| > |z_1
- z_2| > 0$, respectively, and their sums are equal to the branch
\[
f^{p,p}(z_{1}, z_{2})= \sum_{i, j, k, l = 0}^N
a_{ijkl}e^{m_il_p(z_1)}e^{n_jl_p(z_2)}l_p(z_1)^kl_p(z_2)^l(z_1 -
z_2)^{-t}
\]
of $f(z_1, z_2)$ in the region $|z_1| > |z_2| > 0$, the region $|z_2| > |z_1| > 0$, 
the region given by $|z_2| > |z_1- z_2| > 0$ and $|\arg z_{1}-\arg z_{2}|<\frac{\pi}{2}$, respectively.

\item The {\it $L(0)$-grading condition} and {\it $g$-grading condition}: 
Let  $L_{W}^{g}(0)=\res_{x}xY_{W}^g(\omega, x)$.  Then for $n\in \C$ and $\alpha\in \C/\Z$,
$w \in W_{[n]}^{[\alpha]}$,
there exists $K, \Lambda\in \Z_{+}$ such that $(L_{W}^g (0)-n)^{K} w
=(g-e^{2\pi \alpha i})^{\Lambda}w=0$. Moreover,
$gY_{W}^{g}(u,x)v=Y_{W}^{g}(gu,x)gv$.

\item The $L(-1)$-{\it derivative property}: For $v \in V$,
\[
\frac{d}{dx}Y^g_{W}(v, x) = Y^g_{W}(L_{V}(-1)v, x).
\]

\end{enumerate}
A {\it lower-bounded generalized $g$-twisted $V$-module} 
is a $g$-twisted
$V$-module  $W$  such that
for each $n\in \C$, $W_{[n + l]} = 0$ for
sufficiently negative real number $l$. A $g$-twisted $V$-module $W$
is said to be  {\it grading-restricted} if it is lower 
bounded and  for each $n \in
\mathbb{C}$, $\dim W_{[n]}<\infty$.
}
\end{defn}

We shall denote the $g$-twisted
$V$-module just defined by $(W, Y_{W}^{g})$ or simply by $W$ when 
$Y^{g}_{W}$  is clear.

Let $(W, Y^{g}_{W})$ be a  $g$-twisted
$V$-module. Using the notation introduced in Section 2, we have the $p$-th analytic branch $(Y_{W}^{g})^{p}(\omega, z)$ of 
the formal series $Y_{W}^{g}(\omega, x)$ for $p\in \Z$.  Since $g\omega=\omega$, 
$$(Y_{W}^{g})^{p+1}(\omega, z)=(Y_{W}^{g})^{p+1}(g\omega, z)=(Y_{W}^{g})^{p}(\omega, z)$$
for $p\in \Z$. Thus $Y_{W}^{g}(\omega, x)$ involves only integral powers of $x$. 
Let 
$$Y_{W}^{g}(\omega, x)=\sum_{n\in \Z}L_{W}(n)x^{-n-2}.$$
 Then the same argument deriving the Virasoro relations for (untwisted) modules from 
axioms other than those for the Virasoro operators give
$$[L_{W}(m), L_{W}(n)]=(m-n)L_{W}(m+n)+\frac{c}{12}(m^{3}-m)\delta_{m+n, 0}$$
for $m, n\in \Z$, where $c$ is the central charge of $V$. 

Let $(W, Y^{g}_{W})$ be a  $g$-twisted
$V$-module. Let $h$ be an automorphism of $V$ and let
\begin{eqnarray*}
\phi_{h}(Y^{g}): V\times W&\to& W\{x\}[{\rm log} x]\nn
v \otimes w & \mapsto & \phi_{h}(Y^{g})(v, x)w
\end{eqnarray*}
be the linear map defined by
$$ \phi_{h}(Y^{g})(v, x)w=Y^{g}(h^{-1}v, x)w.$$
The following result can be proved by a straightforward use of the axioms:

\begin{prop}
The pair $(W, \phi_{h}(Y^{g}))$ is an $hgh^{-1}$-twisted
$V$-module.
\end{prop}

We shall denote the $hgh^{-1}$-twisted
$V$-module in the proposition above by $\phi_{h}(W)$.

We also need  contragredient twisted $V$-modules.
Let $(W, Y^{g}_{W})$ be a $g$-twisted $V$-module relative to $G$.
Let $W'$ be the graded dual of $W$. Define
a linear map 
\begin{eqnarray*}
(Y_{W}^{g})': V\otimes W' &\to& W'\{x\}[{\rm log} x],\\
v \otimes w' &\mapsto &(Y^g_{W})'(v, x)w'
\end{eqnarray*}
by
$$\langle (Y^g_{W})'(v, x)w', w\rangle=\langle w', Y^{g}_{W}(e^{xL(1)}(-x^{-2})^{L(0)}v, x^{-1})w\rangle$$
for $v\in V$, $w\in W$ and $w'\in W'$. 

\begin{prop}\label{cont-tw-mo}
The pair $(W', (Y^{g}_{W})')$ is a $g^{-1}$-twisted $V$-module.
\end{prop}

The proof of this result is a special case of the proof of Theorem \ref{A} in Section 6 with 
$W_{1}=V$, $g_{1}=1_{V}$, $g_{2}=g$ and $W_{2}=W_{3}=W$ (see also Remark \ref{tw-mod-tw-io}
in Section 4  below). Since the proof 
of Theorem \ref{A} uses only the definition of $(Y^{g}_{W})'$, quoting the proof of  Theorem \ref{A}
to give a proof of Proposition \ref{cont-tw-mo} does not constitute circular reasoning.  

The $g^{-1}$-twisted $V$-module $(W', (Y^{g}_{W})')$ is called the 
{\it contragredient twisted $V$-module of $(W, Y^{g})$}.

\setcounter{equation}{0} \setcounter{thm}{0}

\section{Twisted intertwining operators}

We introduce the notion of  twisted
intertwining operator or intertwining operator among twisted modules in this section. 
Twisted intertwining  operators
in this paper in general involve the logarithm of the variable. Such intertwining 
operators are usually called logarithmic intertwining operators. For simplicity,
we omit the word ``logarithm,''  unless there is a need to emphasize that 
the intertwining operator indeed involve the logarithm of the variable.

\begin{defn}\label{def-tw-int-op}
{\rm Let $g_{1}, g_{2}, g_{3}$ be automorphisms of $V$ and let $W_{1}$, $W_{2}$ and $W_{3}$ be $g_{1}$-, $g_{2}$- 
and $g_{3}$-twisted 
$V$-modules, respectively. A {\it  twisted intertwining operator of type ${W_{3}\choose W_{1}W_{2}}$} is a 
linear map
\begin{eqnarray*}
\mathcal{Y}: W_{1}\otimes W_{2}&\to& W_{3}\{x\}[\log x]\nn
w_{1}\otimes w_{2}&\mapsto& \mathcal{Y}(w_{1}, x)w_{2}=\sum_{k=0}^{K}\sum_{n\in \C}
\mathcal{Y}_{n, k}(w_{1})w_{2}x^{-n-1}(\log x)^{k}
\end{eqnarray*}
satisfying the following conditions:

\begin{enumerate}

\item {\it The lower truncation property}: For $w_{1}\in W_{1}$ and $w_{2}\in W_{2}$, $n\in \C$ and $k=0, \dots, K$, 
$\mathcal{Y}_{n+l, k}(w_{1})w_{2}=0$ for $l\in \N$ sufficiently large.

\item The {\it duality property}:  For $u\in V$, $w_{1}\in W_{1}$, $w_{2}\in W_{2}$
and $w_{3}'\in W_{3}'$, there exists a
multivalued analytic function with preferred branch
\begin{align*}
f&(z_1, z_2; u, w_{1}, w_{2}, w_{3}')\nn
&\quad = \sum_{i, j, k, l, m, n=0}^{N}a_{ijklmn}z_1^{r_{i}}z_2^{s_{j}}(z_1 - z_2)^{t_{k}}(\log z_1)^l
(\log z_2)^{m}(\log (z_1 - z_2))^{n}
\end{align*}
for $N \in \mathbb{N}$, $r_{i}, s_{j}, t_{k}, a_{ijklmn}\in \mathbb{C}$,
such that for $p_{1}, p_{2}, p_{12}\in \Z$, the series 
\begin{align}
\langle w'_{3}, (Y_{W_{3}}
^{g_{3}})^{p_{1}}(u, z_{1})\mathcal{Y}^{p_{2}}(w_{1}, z_{2})w_{2}\rangle
&=\sum_{n\in \C}\langle w'_{3}, (Y_{W_{3}}
^{g_{3}})^{p_{1}}(u, z_{1})\pi_{n}
\mathcal{Y}^{p_{2}}(w_{1}, z_{2})w_{2}\rangle,
\label{int-prod}\\
\langle w'_{3}, \mathcal{Y}^{p_{2}}(w_{1}, z_{2})
(Y_{W_{2}}^{g_{2}})^{p_{1}}(u, z_{1})w_{2}\rangle
&=\sum_{n\in \C}\langle w'_{3}, \mathcal{Y}^{p_{2}}(w_{1}, z_{2})
\pi_{n}(Y_{W_{2}}^{g_{2}})^{p_{1}}(u, z_{1})w_{2}\rangle,
\label{int-rev-prod}\\
 \langle w'_{3}, \mathcal{Y}^{p_{2}}(
(Y_{W_{1}}^{g_{1}})^{p_{12}}(u, z_{1}-z_{2})v, z_{2})w\rangle
&=\sum_{n\in \C}\langle w'_{3}, \mathcal{Y}^{ p_{2}}(\pi_{n}
(Y_{W_{1}}^{g_{1}})^{p_{12}}(u, z_{1}-z_{2})w_{1}, z_{2})w_{2}\rangle\label{int-iter}
\end{align}
are absolutely convergent in the regions 
$|z_{1}|>|z_{2}|>0$,
$|z_{2}|>|z_{1}|>0$, $|z_{2}|>|z_{1}-z_{2}|>0$, respectively. Moreover, their sums are equal 
to the branches
\begin{align*}
f&^{p_{1}, p_{2}, p_{1}}(z_{1}, z_{2}; u, w_{1}, w_{2}, w_{3}')\nn
&=\sum_{i, j, k, l, m, n=0}^{N}a_{ijklmn}e^{r_{i}l_{p_{1}}(z_{1})}
e^{s_{j}l_{p_{2}}(z_2)}e^{t_{k}l_{p_{1}}(z_1 - z_2)}(l_{p_{1}}(z_1))^l
(l_{p_{2}}(z_2))^{m}(l_{p_{1}}(z_1 - z_2))^{n},\nn
f&^{p_{1}, p_{2}, p_{2}}(z_{1}, z_{2}; u, w_{1}, w_{2}, w_{3}')\nn
&=\sum_{i, j, k, l, m, n=0}^{N}a_{ijklmn}e^{r_{i}l_{p_{1}}(z_{1})}
e^{s_{j}l_{p_{2}}(z_2)}e^{t_{k}l_{p_{2}}(z_1 - z_2)}(l_{p_{1}}(z_1))^l
(l_{p_{2}}(z_2))^{m}(l_{p_{2}}(z_1 - z_2))^{n},\nn
f&^{p_{2}, p_{2}, p_{12}}(z_{1}, z_{2}; u, w_{1}, w_{2}, w_{3}')\nn
&=\sum_{i, j, k, l, m, n=0}^{N}a_{ijklmn}e^{r_{i}l_{p_{2}}(z_{1})}
e^{s_{j}l_{p_{2}}(z_2)}e^{t_{k}l_{p_{12}}(z_1 - z_2)}(l_{p_{2}}(z_1))^l
(l_{p_{2}}(z_2))^{m}(l_{p_{12}}(z_1 - z_2))^{n},
\end{align*}
respectively, of $f(z_1, z_2; u, w_{1}, w_{2}, w_{3}')$ (recall the notations and convention in Section 2)  in the region given by 
$|z_{1}|>|z_{2}|>0$ and $|\arg (z_{1}-z_{2})-\arg z_{1}|<\frac{\pi}{2}$,
the region given by $|z_{2}|>|z_{1}|>0$ and $-\frac{3\pi}{2}< \arg (z_{1}-z_{2})-\arg z_{2}<-\frac{\pi}{2}$,
the region given by $|z_{2}|>|z_{1}-z_{2}|>0$ and $|\arg  z_{1}- \arg z_{2}|<\frac{\pi}{2}$, respectively.

\item The {\it $L(-1)$-derivative property}: 
$$\frac{d}{dx}\mathcal{Y}(w_{1}, x)=\mathcal{Y}(L(-1)w_{1}, x).$$
\end{enumerate}
 }
\end{defn}

\begin{rema}\label{tw-mod-tw-io}
{\rm Let $(W, Y_{W}^{g})$ be a $g$-twisted $V$-module. Then by definition, 
the vertex operator map $Y_{W}^{g}$
is a twisted intertwining operator of type ${W\choose VW}$.}
\end{rema}

\begin{rema}\label{int-op-for-sub-m}
{\rm Let $g_{1}$, $g_{2}$ and $g_{3}$  be automorphisms of $V$, let $V^{\langle g_{1}, g_{2}, g_{3}\rangle}$ be the 
vertex operator subalgebra of $V$ consisting of elements of $V$ fixed under the actions of $g_{1}, g_{2}$ 
and $g_{3}$. 
Let $W_{1}$, $W_{2}$ and $W_{3}$ be $g_{1}$-, $g_{2}$- 
and $g_{3}$-twisted 
$V$-modules and let $\mathcal{Y}$ be a twisted intertwining operator of type ${W_{3}\choose W_{1}W_{2}}$. 
Then $W_{1}$, $W_{2}$ and $W_{3}$  are $V^{\langle g_{1}, g_{2}, g_{3}\rangle}$-modules and 
$\mathcal{Y}$ is an (untwisted)  intertwining operator of type ${W_{3}\choose W_{1}W_{2}}$ when 
$W_{1}$, $W_{2}$ and $W_{3}$  are viewed as $V^{\langle g_{1}, g_{2}, g_{3}\rangle}$-modules.
But in general, an (untwisted)  intertwining operator of type ${W_{3}\choose W_{1}W_{2}}$ when 
$W_{1}$, $W_{2}$ and $W_{3}$  are viewed as $V^{\langle g_{1}, g_{2}, g_{3}\rangle}$-modules
might not be a twisted intertwining operator of type ${W_{3}\choose W_{1}W_{2}}$. In fact, an (untwisted)  intertwining operator of type ${W_{3}\choose W_{1}W_{2}}$ when 
$W_{1}$, $W_{2}$ and $W_{3}$  are viewed as $V^{\langle g_{1}, g_{2}, g_{3}\rangle}$-modules 
is required to satisfy only the duality property for vertex operators associated to elements in $V^{\langle g_{1}, g_{2}, g_{3}\rangle}$
while a twisted intertwining operator of type ${W_{3}\choose W_{1}W_{2}}$
must satisfy the more restrictive duality property in Definition \ref{def-tw-int-op}. 
This is the reason why even when $W_{1}$, $W_{2}$ and $W_{3}$ are known to be twisted $V$-modules, we still 
want to add the word ``twisted" in front of ``intertwining operator" to call an intertwining 
operator of type ${W_{3}\choose W_{1}W_{2}}$ in Definition \ref{def-tw-int-op}
a twisted intertwining operator.}
\end{rema}

\begin{rema}
{\rm In the duality property in the definition above, we require that the sum of  (\ref{int-rev-prod}) 
is equal to $f^{p_{1}, p_{2}, p_{2}}(z_{1}, z_{2}; u, w_{1}, w_{2}, w_{3}')$ in the region given by
$|z_{2}|>|z_{1}|>0$ and $-\frac{3\pi}{2}< \arg (z_{1}-z_{2})-\arg z_{2}<-\frac{\pi}{2}$. 
The choice of this region in fact gives an order of $W_{1}$ and $W_{2}$ to be $W_{1}$ first and $W_{2}$ second. 
See Theorem \ref{g3=g1g2} and especially its proof below.
The other region is the region given by 
$|z_{2}|>|z_{1}|>0$ and $\frac{\pi}{2}< \arg (z_{1}-z_{2})-\arg z_{2}<\frac{3\pi}{2}$. If we require 
the sum of  (\ref{int-rev-prod}) 
is equal to $f^{p_{1}, p_{2}, p_{2}}(z_{1}, z_{2}; u, w_{1}, w_{2}, w_{3}')$ in this region, then 
the order of $W_{1}$ and $W_{2}$ is chosen to be $W_{2}$ first and $W_{1}$ second.
We choose the more natural order. }
\end{rema}

We shall need the following two lemmas:

\begin{lemma}\label{int-rev-prod-1}
In Definition \ref{def-tw-int-op}, the requirement in the duality property
that the sum of (\ref{int-rev-prod})
be equal to $f^{p_{1}, p_{2}, p_{2}}(z_{1}, z_{2}; u, w_{1}, w_{2}, w_{3}')$ in the region given
by $|z_{2}|>|z_{1}|>0$ and $-\frac{3\pi}{2}< \arg (z_{1}-z_{2})-\arg z_{2}<-\frac{\pi}{2}$
can be replaced by the requirement that the sum of (\ref{int-rev-prod})
be equal to $f^{p_{1}, p_{2}, p_{2}-1}(z_{1}, z_{2}; u, w_{1}, w_{2}, w_{3}')$ in the region given
by $|z_{2}|>|z_{1}|>0$ and $\frac{\pi}{2}< \arg (z_{1}-z_{2})-\arg z_{2}<\frac{3\pi}{2}$. 
In the same definition,  the requirement in the duality property that the sum of (\ref{int-prod})
be equal to $f^{p_{1}, p_{2}, p_{1}}(z_{1}, z_{2}; u, w_{1}, w_{2}, w_{3}')$ in the region given
by $|z_{1}|>|z_{2}|>0$ and $|\arg (z_{1}-z_{2})-\arg z_{2}|<\frac{\pi}{2}$
can be replaced by the requirement that the sum of (\ref{int-prod})
be equal to $f^{p_{1}, p_{2}, p_{1}-1}(z_{1}, z_{2}; u, w_{1}, w_{2}, w_{3}')$ in the region given
by $|z_{2}|>|z_{1}|>0$ and $-2\pi< \arg (z_{1}-z_{2})-\arg z_{2}<-\frac{3\pi}{2}$ and to 
$f^{p_{1}, p_{2}, p_{1}+1}(z_{1}, z_{2}; u, w_{1}, w_{2}, w_{3}')$ in the region given
by $|z_{2}|>|z_{1}|>0$ and $\frac{3\pi}{2}< \arg (z_{1}-z_{2})-\arg z_{2}<2\pi$. 
\end{lemma}
\pf
We prove only the first part. The second part can be proved similarly. 

\begin{figure}[h]
\centering
\includegraphics{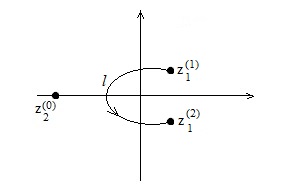}
\caption{The path $l$\label{figure0}}
\end{figure}

Assume that $\Y$ is a twisted intertwining operator satisfying Definition \ref{def-tw-int-op}.
We choose a path $l$ in the region given by $|z_{1}|>|z_{2}|>0$
from a point $(z^{(1)}_{1}, z^{(0)}_{2})$ in the subregion given by 
$|z_{2}|>|z_{1}|>0$ and $-\frac{3\pi}{2}< \arg (z_{1}-z_{2})-\arg z_{2}<-\frac{\pi}{2}$
to a point $(z_{1}^{(2)}, z_{2}^{(0)})$ in the subregion given by 
$|z_{2}|>|z_{1}|>0$ and $\frac{\pi}{2}< \arg (z_{1}-z_{2})-\arg z_{2}<\frac{3\pi}{2}$
by letting $z_{1}$ pass through the set given by $\arg (z_{1}-z_{2})=0$ in the counter clockwise
direction for the variable $z_{1}-z_{2}$ but keeping $\arg z_{1}$ between $0$ and $2 \pi$ and
fixing $z_{2}=z_{2}^{(0)}$. See Figure \ref{figure0}.
The sum of the series (\ref{int-rev-prod}) is an analytic function of $z_{1}$ and $z_{2}$ and thus we can analytically 
extend its value at $(z_{1}^{(1)}, z_{2}^{(0)})$ through the path $l$ to its value at $(z_{1}^{(2)}, z_{2}^{(0)})$.
At $(z_{1}^{(1)}, z_{2}^{(0)})$, its value is given by $f^{p_{1}, p_{2}, p_{2}}(z_{1}^{(1)}, z_{2}^{(0)}; 
u, w_{1}, w_{2}, w_{3}')$. When the path $l$ pass the point at which $\arg (z_{1}-z_{2})=0$, there is a 
jump of $\arg (z_{1}-z_{2})$ from $0$ to $2\pi$. When $\arg (z_{1}-z_{2})=0$, its value (at $z_{1},z_{2})$  is still 
$f^{p_{1}, p_{2}, p_{2}}(z_{1}, z_{2}; 
u, w_{1}, w_{2}, w_{3}')$. But after the jump, since the sum is analytic and in particular is continuous, 
its value at $(z_{1}, z_{2})$  must be $f^{p_{1}, p_{2}, p_{2}-1}(z_{1}, z_{2}; u, w_{1}, w_{2}, w_{3}')$. In particular, 
its value at the arbitrary point $(z_{1}^{(2)}, z_{2}^{(0)})$ in the region given by 
$|z_{2}|>|z_{1}|>0$ and $\frac{\pi}{2}< \arg (z_{1}-z_{2})-\arg z_{2}<\frac{3\pi}{2}$ must be 
$f^{p_{1}, p_{2}, p_{2}-1}(z_{1}^{(2)}, z_{2}^{(0)}; u, w_{1}, w_{2}, w_{3}')$. 

If we assume that $\Y$ satisfies the
requirement that the sum of (\ref{int-rev-prod})
be equal to 
$$f^{p_{1}, p_{2}, p_{2}-1}(z_{1}, z_{2}; u, w_{1}, w_{2}, w_{3}')$$
in the region given
by $|z_{2}|>|z_{1}|>0$ and $\frac{\pi}{2}< \arg (z_{1}-z_{2})-\arg z_{2}<\frac{3\pi}{2}$ and all
the other axioms in Definition \ref{def-tw-int-op}, a completely analogous argument shows 
that  the sum of (\ref{int-rev-prod})
be equal to $f^{p_{1}, p_{2}, p_{2}}(z_{1}, z_{2}; u, w_{1}, w_{2}, w_{3}')$ in the region given
by $|z_{2}|>|z_{1}|>0$ and $-\frac{3\pi}{2}< \arg (z_{1}-z_{2})-\arg z_{2}<-\frac{\pi}{2}$.
\epfv

\begin{lemma}
For $p_{1}, p_{2}, p_{12}\in \Z$, $u\in V$, $w_{1}\in W_{1}$, $w_{2}\in W_{2}$
and $w_{3}'\in W'_{3}$, we have 
\begin{equation}\label{skew-sym-8.5}
f^{p_{1}, p_{2}, p_{12}+1}(z_1, z_2; g_{1}u, w_{1}, w_{2},w_{3}')
=f^{p_{1}, p_{2}, p_{12}}(z_1-z_{2}, z_2; u, w_{1}, w_{2}, w_{3}')
\end{equation}
and 
\begin{equation}\label{skew-sym-8-1.5}
f^{p_{1}+1, p_{2}, p_{12}}(z_1, z_2; g_{2}u, w_{1}, w_{2}, w_{3}')
=f^{p_{1}, p_{2}, p_{12}}(z_1, z_2; u, w_{1}, w_{2}, w_{3}'),
\end{equation}
where $f^{p_{1}, p_{2}, p_{12}}(z_{1}, z_{2}; u, w_{1}, w_{2}, w_{3}')$
for $p_{1}, p_{2}, p_{12}\in \Z$ is the branch given by  $p_{1}, p_{2}, p_{12}$
of the multivalued analytic function $f(z_{1}, z_{2}; u, w_{1}, w_{2}, w_{3}')$  
with preferred branch in
Definition \ref{def-tw-int-op}.
\end{lemma}
\pf
By the duality property for $\Y$, when $|z_{2}|>|z_{1}-z_{2}|>0$ and 
$|\arg z_{1}-\arg z_{2}|<\frac{\pi}{2}$,
$$\langle w'_{3}, \mathcal{Y}^{p_{2}}((Y_{W_{1}}
^{g_{1}})^{p_{12}+1}(g_{1}u, z_{1}-z_{2})w_{1}, z_{2})
w_{2}\rangle$$
converges absolutely to 
$f^{p_{2}, p_{2}, p_{12}+1}(z_1, z_2; g_{1}u, w_{1}, w_{2}, w_{3}')$
and 
$$\langle w'_{3}, \mathcal{Y}^{p_{2}}((Y_{W_{1}}
^{g_{1}})^{p_{12}}(u, z_{1}-z_{2})w_{1}, z_{2})
w_{2}\rangle$$
converges absolutely to 
$f^{p_{2}, p_{2}, p_{12}}(z_1, z_2; u, w_{1}, w_{2}, w_{3}')$.
But
$$\langle w'_{3}, \mathcal{Y}^{p_{2}}((Y_{W_{1}}
^{g_{1}})^{p_{12}+1}(g_{1}u, z_{1}-z_{2})w_{1}, z_{2})
w_{2}\rangle
=\langle w'_{3}, \mathcal{Y}^{p_{2}}((Y_{W_{1}}
^{g_{1}})^{p_{12}}(u, z_{1}-z_{2})w_{1}, z_{2})
w_{2}\rangle$$
Thus we have 
$$f^{p_{2}, p_{2}, p_{12}+1}(z_1, z_2; g_{1}u, w_{1}, w_{2}, w_{3}')
=f^{p_{2}, p_{2}, p_{12}}(z_1, z_2; u, w_{1}, w_{2}, w_{3}').$$
For general $p_{1}, p_{2}, p_{12}\in \Z$,  we obtain (\ref{skew-sym-8.5}) by analytic extensions.

On the other hand, by the duality property for $\Y$, when $|z_{2}|>|z_{1}|>0$ and 
$-\frac{3\pi}{2}<\arg z_{1}-\arg (-z_{2})<-\frac{\pi}{2}$,
$$\langle w'_{3}, \mathcal{Y}^{p_{2}}(w_{1}, z_{2})
(Y_{W_{1}}^{g_{2}})^{p_{1}+1}(g_{2}u, z_{1})w_{2}\rangle$$
converges absolutely to 
$f^{p_{1}+1, p_{2}, p_{2}}(z_1, z_2; g_{2}u, w_{1}, w_{2}, w_{3}')$
and 
$$\langle w'_{3}, \mathcal{Y}^{p_{2}}(w_{1}, z_{2})
(Y_{W_{1}}^{g_{2}})^{p_{1}}(u, z_{1})w_{2}\rangle$$
converges absolutely to 
$f^{p_{1}, p_{2}, p_{2}}(z_1, z_2; u, w_{1}, w_{2}, w_{3}')$.
But
$$\langle w'_{3}, \mathcal{Y}^{p_{2}}(w_{1}, z_{2})
(Y_{W_{1}}^{g_{2}})^{p_{1}+1}(g_{2}u, z_{1})w_{2}\rangle
=\langle w'_{3}, \mathcal{Y}^{p_{2}}(w_{1}, z_{2})
(Y_{W_{1}}^{g_{2}})^{p_{1}}(u, z_{1})w_{2}\rangle.$$
Thus we have 
$$f^{p_{1}+1, p_{2}, p_{2}}(z_1, z_2; g_{2}u, w_{1}, w_{2}, w_{3}')
=f^{p_{1}, p_{2}, p_{2}}(z_1-z_{2}, -z_2; u, w_{1}, w_{2}, w_{3}').$$
For general $p_{1}, p_{2}, p_{12}\in \Z$,  we obtain (\ref{skew-sym-8-1.5})  by analytic extensions.
\epfv

We now prove that under suitable minor conditions, $g_{3}=g_{1}g_{2}$ for the twisted intertwining 
operator defined in Definition  \ref{def-tw-int-op}. 

\begin{figure}[h]
\centering
\includegraphics{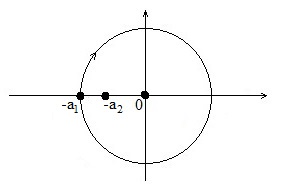}
\caption{The loop $\Gamma_{1}$\protect\label{figure1}}
\end{figure}

\begin{figure}[h]
\centering
\includegraphics{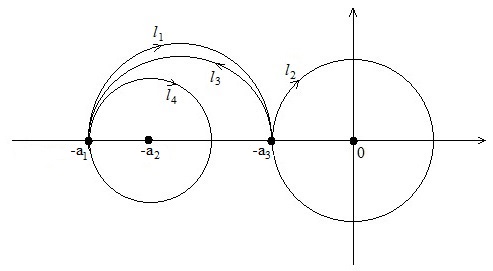}
\caption{The loop $\Gamma_{2}$\label{figure2}}
\end{figure}

\begin{thm}\label{g3=g1g2}
Let $g_{1}, g_{2}, g_{3}$ be automorphisms of $V$ and let $W_{1}$, $W_{2}$ and $W_{3}$ be $g_{1}$-, $g_{2}$- 
and $g_{3}$-twisted 
$V$-modules, respectively.
Assume that the vertex operator map for $W_{3}$ given by $u\mapsto Y_{W_{3}}^{g_{3}}(u, x)$
is injective. If there exists a twisted intertwining operator $\mathcal{Y}$ of type ${W_{3}\choose W_{1}W_{2}}$ such that 
the coefficients of the series $\mathcal{Y}(w_{1}, x)w_{2}$ for $w_{1}\in W_{1}$ and $w_{2}\in W_{2}$ span
$W_{3}$,
then $g_{3}=g_{1}g_{2}$.
\end{thm}
\pf
Let $\mathcal{Y}$ be a twisted  intertwining operator of type ${W_{3}\choose W_{1}W_{2}}$  such that 
the coefficients of $\mathcal{Y}(w_{1}, x)w_{2}$ for $w_{1}\in W_{1}$ and $w_{2}\in W_{2}$ span
$W_{3}$.
For $u\in V$, $w_{1}\in W_{1}$, $w_{2}\in W_{2}$
and $w_{3}'\in W_{3}'$, consider the multivalued analytic function $f(z_{1}, z_{2}; u, w_{1}, w_{2}, w_{3}')$ 
with preferred branch for the twisted intertwining operator $\mathcal{Y}$ (see Definition \ref{def-tw-int-op}). 
Fix $z_{2}$ to be a nonzero negative real number $-a_{2}$ where 
$a_{2}\in \R_{+}$.
Then for any $p\in \Z$, we have an  analytic function 
$$f_{-a_{2}}^{p}(z_{1})=\sum_{i, j, k, l, m, n=0}^{N}a_{ijklmn}e^{r_{i}l_{p}(z_{1})}
e^{s_{j}l_{p}(-a_2)}e^{t_{k}l_{p}(z_1 +a_{2})}(l_{p}(z_1))^l
(l_{p}(-a_2))^{m}(l_{p}(z_1 +a_2))^{n}$$
of $z_{1}$ and can be analytically extended to a multivalued analytic function
$$f_{- a_{2}}(z_{1})=\sum_{i, j, k, l, m, n=0}^{N}a_{ijklmn}z_1^{r_i}e^{s_{j}l_{p}(-a_2)}(z_1 + a_2)^{t_{k}}(\log z_1)^l
(l_{p}(-a_2))^{m}(\log (z_1 +a_2))^{n}$$
of $z_{1}$ with preferred branch.
Let $a_{1}\in \R_{+}$  such that $a_{1}>a_{2}>a_{1}-a_{2}$.
Consider the loop $\Gamma_{1}$  in the $z_{1}$ plane based at $z_{1}=-a_{1}$ 
in Figure \ref{figure1}.
\noindent We consider the  value $f^{p}_{-a_{2}}(-a_{1})$ of the multivalued analytic function $f_{-a_{2}}(z_{1})$ at $-a_{1}$.
By the definition of twisted intertwining operator above, 
\begin{equation}\label{g3=g1g2-1}
f^{p}_{-a_{2}}(-a_{1})=\langle w'_{3}, (Y_{W_{3}}
^{g_{3}})^{p}(u, -a_{1})\mathcal{Y}^{p}(w_{1}, -a_{2})w_{2}\rangle.
\end{equation}
But by definition, when $z_{1}$ goes around the loop above, the right-hand side of 
(\ref{g3=g1g2-1}) changes to 
\begin{equation}\label{g3=g1g2-2}
\langle w'_{3}, (Y_{W_{3}}
^{g_{3}})^{p-1}(u, -a_{1})\mathcal{Y}^{p}(w_{1}, -a_{2})w_{2}\rangle.
\end{equation}
By the equivariance property of the $g_{3}$-twisted module $W_{3}$,
(\ref{g3=g1g2-2}) is equal to 
\begin{equation}\label{g3=g1g2-3}
\langle w'_{3}, (Y_{W_{3}}
^{g_{3}})^{p}(g_{3}u, -a_{1})\mathcal{Y}^{p}(w_{1}, -a_{2})w_{2}\rangle.
\end{equation}
We also consider another loop $\Gamma_{2}$  in the $z_{1}$ plane and based at $z_{1}=-a_{1}$ 
given in the order $l_{1}$ first, $l_{2}$ second, $l_{3}$ third and $l_{4}$ last
in Figure \ref{figure2}.
The loop $\Gamma_{2}$ is in fact homotopy equivalent to the loop $\Gamma_{1}$. Thus when 
$z_{1}$ goes around $\Gamma_{2}$, the right-hand side of 
(\ref{g3=g1g2-1}) also changes to (\ref{g3=g1g2-3}). 

On the other hand, we look at how the function values change when $z_{1}$ goes through
$l_{1}$, $l_{2}$, $l_{3}$ and $l_{4}$. When $z_{1}$ goes from $-a_{1}$ to $-a_{3}$ (see Figure \ref{figure2})
through 
$l_{1}$,  the right-hand side of 
(\ref{g3=g1g2-1})  changes to 
\begin{equation}\label{g3=g1g2-4}
\sum_{i, j, k, l, m, n=0}^{N}a_{ijklmn}e^{r_{i}l_{p}(-a_{3})}
e^{s_{j}l_{p}(-a_2)}e^{t_{k}l_{p}(-a_{3} +a_{2})}(l_{p}(-a_{3}))^l
(l_{p}(-a_2))^{m}(l_{p}(-a_{3} +a_2))^{n}
\end{equation}
Note that
$\arg (-a_{2})=\pi$, $\arg (-a_{3})=\pi$ and $\arg (-a_{3}+a_{2})=0$. Hence
$\arg (-a_{3}+a_{2})-\arg (-a_{2})=-\pi$. Since $|-a_{1}|>|-a_{3}|>0$ and 
$-\frac{3\pi}{2}<\arg (-a_{3}+a_{2})-\arg (-a_{2})<-\frac{\pi}{2}$,   by the duality property of 
the twisted intertwining operator $\Y$, (\ref{g3=g1g2-4}) is equal to 
\begin{equation}\label{g3=g1g2-5}
\langle w'_{3}, \mathcal{Y}^{p}(w_{1}, -a_{2})(Y_{W_{2}}
^{g_{2}})^{p}(u, -a_{3})w_{2}\rangle.
\end{equation}
Next let $z_{1}$ go around the loop $l_{2}$. Then (\ref{g3=g1g2-5}) changes to 
\begin{equation}\label{g3=g1g2-6}
\langle w'_{3}, \mathcal{Y}^{p}(w_{1}, -a_{2})(Y_{W_{2}}
^{g_{2}})^{p-1}(u, -a_{3})w_{2}\rangle.
\end{equation}
By the equivariance property of the $g_{2}$-twisted module $W_{2}$,
(\ref{g3=g1g2-6}) is equal to 
\begin{equation}\label{g3=g1g2-7}
\langle w'_{3}, \mathcal{Y}^{p}(w_{1}, -a_{2})(Y_{W_{2}}
^{g_{2}})^{p}(g_{2}u, -a_{3})w_{2}\rangle.
\end{equation}
Now let $z_{1}$ go from $-a_{3}$ to $-a_{1}$ through $l_{3}$. 
Then by reversing the argument above on the change of the values 
when $z_{1}$ goes through $l_{1}$, we see that the values change from 
(\ref{g3=g1g2-7}) to
\begin{equation}\label{g3=g1g2-8}
\langle w'_{3}, (Y_{W_{3}}
^{g_{3}})^{p}(g_{2}u, -a_{1})\mathcal{Y}^{p}(w_{1}, -a_{2})w_{2}\rangle.
\end{equation}
Since $|-a_{1}|>|-a_{2}|>|-a_{1}-(-a_{2})|>0$, $|\arg (-a_{1}-(a_{2}))-\arg (-a_{1})|=0<\frac{\pi}{2}$ 
 and $|\arg (-a_{1})-\arg (-a_{2})|=0<\frac{\pi}{2}$,
by the duality property of the twisted intertwining operator $\Y$, (\ref{g3=g1g2-8})
is equal to 
\begin{equation}\label{g3=g1g2-9}
\langle w'_{3}, \mathcal{Y}^{p}((Y_{W_{1}}
^{g_{3}})^{p}(g_{2}u, -a_{1}+a_{2})w_{1}, -a_{2})w_{2}\rangle.
\end{equation}
Finally let $z_{1}$ go around the loop $l_{4}$. The value changes from 
(\ref{g3=g1g2-9}) to 
\begin{equation}\label{g3=g1g2-10}
\langle w'_{3}, \mathcal{Y}^{p}((Y_{W_{1}}
^{g_{3}})^{p-1}(g_{2}u, -a_{1}+a_{2})w_{1}, -a_{2})w_{2}\rangle.
\end{equation}
By the equivariance property of the $g_{1}$-twisted module $W_{1}$,
(\ref{g3=g1g2-10}) is equal to 
\begin{equation}\label{g3=g1g2-11}
\langle w'_{3}, \mathcal{Y}^{p}((Y_{W_{1}}
^{g_{3}})^{p}(g_{1}g_{2}u, -a_{1}+a_{2})w_{1}, -a_{2})w_{2}\rangle.
\end{equation}
Again since $|-a_{1}|>|-a_{2}|>|-a_{1}-(-a_{2})|>0$, $|\arg (-a_{1}-(a_{2}))-\arg (-a_{1})|=0<\frac{\pi}{2}$ 
 and $|\arg (-a_{1})-\arg (-a_{2})|=0<\frac{\pi}{2}$,
by the duality property of the twisted intertwining operator $\Y$, (\ref{g3=g1g2-11})
is equal to 
\begin{equation}\label{g3=g1g2-12}
\langle w'_{3}, (Y_{W_{3}}
^{g_{3}})^{p}(g_{1}g_{2}u, -a_{1})\mathcal{Y}^{p}(w_{1}, -a_{2})w_{2}\rangle.
\end{equation}

From the discussions above, we see that when 
$z_{1}$ goes around $\Gamma_{2}$, the right-hand side of 
(\ref{g3=g1g2-1}) changes to 
(\ref{g3=g1g2-12}). Thus we see that 
(\ref{g3=g1g2-3}) and (\ref{g3=g1g2-12}) are  equal, that is
\begin{equation}\label{g3=g1g2-13}
\langle w'_{3}, (Y_{W_{3}}
^{g_{3}})^{p}(g_{3}u-g_{1}g_{2}u, -a_{1})\mathcal{Y}^{p}(w_{1}, -a_{2})w_{2}\rangle=0.
\end{equation}
Since $w_{1}$, $w_{2}$ and $w'_{3}$ are arbitrary and the coefficients of 
$\mathcal{Y}(w_{1}, x)w_{2}$ for $w_{1}\in W_{1}$ and $w_{2}\in W_{2}$ span
$W_{3}$, we obtain from (\ref{g3=g1g2-13})
\begin{equation}\label{g3=g1g2-14}
(Y_{W_{3}}
^{g_{3}})^{p}(g_{3}u-g_{1}g_{2}u, -a_{1})=0.
\end{equation}
Replacing $u$ in  (\ref{g3=g1g2-14}) by $L_{V}(-1)^{m}u$ for $m\in \N$, using the fact that $L_{V}(-1)$ commutes with $g_{1}$, $g_{2}$
and $g_{3}$ and then using
the $L(-1)$-derivative property for the twisted vertex operator $Y^{g_{3}}_{W_{3}}$, we
obtain
$$
\frac{d^{m}}{dx^{m}}Y_{W_{3}}
^{g_{3}}(g_{3}u-g_{1}g_{2}u, x)\lbar_{x^{n}=e^{nl_{p}(-a_{1})},\; \log x=l_{p}(-a_{1})}=0.
$$
Using the Taylor series expansion, we obtain
\begin{equation}\label{g3=g1g2-16}
Y_{W_{3}}
^{g_{3}}(g_{3}u-g_{1}g_{2}u, x)\lbar_{x^{n}=e^{nl_{p}(z)},\; \log x=l_{p}(z)}=0
\end{equation}
for $z\in \C^{\times}$. Thus
\begin{equation}\label{g3=g1g2-15}
Y_{W_{3}}
^{g_{3}}(g_{3}u-g_{1}g_{2}u, x)=0.
\end{equation}
Since the vertex operator map $u\mapsto Y_{W_{3}}
^{g_{3}}(u, x)$ is injective,  (\ref{g3=g1g2-15}) implies 
$g_{3}u-g_{1}g_{2}u=0$ or $g_{3}u=g_{1}g_{2}u$. 
Since $u$ is also arbitrary, we obtain $g_{3}=g_{1}g_{2}$. 
\epfv

\begin{rema}
{\rm The proof of Theorem \ref{g3=g1g2} can be intuitively
understood using two braiding graphs. See Figure \ref{figure3}.
\begin{figure}[h]
\centering
\includegraphics{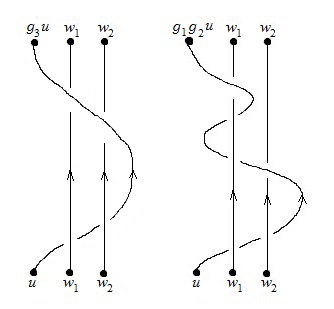}
\caption{The braiding graphs corresponding to $\Gamma_{1}$ (left) and 
$\Gamma_{2}$ (right)\label{figure3}}
\end{figure}
Since the two braiding graphs are topologically equivalent (isotopic),
the corresponding algebraic objects are equal and thus we have $g_{3}u=g_{1}g_{2}u$.
In fact, just like the theory of braided tensor categories, the correspondence 
between algebraic and analytic calculations and the braiding graphs can be 
made mathematically precise so that proofs such as the one for Theorem \ref{g3=g1g2}
can be given using such graphs.  Note that in these graphs, we have suppressed 
the associativity for the twisted vertex operators and the twisted intertwining operator, just as what 
people usually do in the graphs for braided tensor categories. We also note that these braiding graphs 
explain only those results that are topological in nature. For analytic results such as those we 
shall prove in the next two sections,
these graphs are not very useful.
}
\end{rema}

Because of Theorem \ref{g3=g1g2}, in the rest of this paper, we shall discuss only 
twisted intertwining operators of type ${W_{3}\choose W_{1}W_{2}}$ with $W_{1}$, 
$W_{2}$ and $W_{3}$ being $g_{1}$-, $g_{2}$- 
and $g_{1}g_{2}$-twisted 
$V$-modules.

\setcounter{equation}{0} \setcounter{thm}{0}

\section{The skew symmetry isomorphisms}

In this section, we construct what we call the skew-symmetry isomorphisms between 
the spaces of twisted intertwining operators of suitable types. These linear isomorphisms
correspond to braidings in the still-to-be-constructed $G$-crossed braided
tensor category structure on the category of $g$-twisted $V$-modules for all $g$ in a 
group $G$ of automorphisms of $V$.

Let $g_{1}, g_{2}$ be automorphisms of $V$,  $W_{1}$, $W_{2}$ and $W_{3}$ $g_{1}$-, $g_{2}$- 
and $g_{1}g_{2}$-twisted 
$V$-modules and $\mathcal{Y}$ a twisted intertwining operator 
of type ${W_{3}\choose W_{1}W_{2}}$. 
We define  linear maps
\begin{eqnarray*}
\Omega_{\pm}(\Y): W_{2}\otimes W_{1}&\to& W_{3}\{x\}[\log x]\nn
w_{2}\otimes w_{1}&\mapsto& \Omega_{\pm}(\Y)(w_{2}, x)w_{1}
\end{eqnarray*}
by
\begin{equation}\label{omega1}
\Omega_{\pm}(\Y)(w_{2}, x)w_{1}=e^{xL(-1)}\Y(w_{1}, y)w_{2}\lbar_{y^{n}=e^{\pm n\pi \i}x^{n}, \;
\log y=\log x\pm \pi \i}
\end{equation}
for $w_{1}\in W_{1}$ and $w_{2}\in W_{2}$. 

From the definition (\ref{omega1}), for $p\in \Z$, $w_{1}\in W_{1}$, $w_{2}\in W_{2}$ and $z\in \C^{\times}$, 
\begin{align*}
\Omega_{\pm}(\Y)^{p}(w_{2}, z)w_{1}&=\Omega_{\pm}(\Y)(w_{2}, x)w_{1}\lbar_{x^{n}=e^{nl_{p}(z)},\;\log x=l_{p}(z)}\nn
&=\left(e^{xL(-1)}\Y(w_{1}, y)w_{2}\lbar_{y^{n}=e^{\pm n\pi \i}x^{n}, \;
\log y=\log x\pm \pi \i}\right)\lbar_{x^{n}=e^{nl_{p}(z)},\;\log x=l_{p}(z)}\nn
&=e^{zL(-1)}\Y(w_{1}, y)w_{2}\lbar_{y^{n}=e^{n(l_{p}(z)\pm \pi \i)}, \;\log y=l_{p}(z)\pm \pi \i}.
\end{align*}
When $\arg z<\pi$ and $\arg z\ge \pi$, $\arg (-z)=\arg z+\pi$ and $\arg (-z)=\arg z-\pi$, respectively. 
Hence 
$$e^{zL(-1)}\Y(w_{1}, y)w_{2}\lbar_{y^{n}=e^{n(l_{p}(z)+ \pi \i)}, \;\log y=l_{p}(z)+ \pi \i}=
e^{zL(-1)}\Y^{p}(w_{1}, -z)w_{2}$$  
when $\arg z<\pi$ and 
$$e^{zL(-1)}\Y(w_{1}, y)w_{2}\lbar_{y^{n}=e^{n(l_{p}(z)- \pi \i)}, \;\log y=l_{p}(z)- \pi \i}=
e^{zL(-1)}\Y^{p}(w_{1}, -z)w_{2}$$  
when $\arg z\ge \pi$.  In particular, for $w_{1}\in W_{1}$, $w_{2}\in W_{2}$ and $z\in \C^{\times}$ satisfying
$\arg z<\pi$ and $\arg z \ge \pi$, we have
\begin{equation}\label{omega2}
\Omega_{+}(\Y)^{p}(w_{2}, z)w_{1}=e^{zL(-1)}\Y^{p}(w_{1}, -z)w_{2}.
\end{equation}
and 
\begin{equation}\label{omega3}
\Omega_{-}(\Y)^{p}(w_{2}, z)w_{1}=e^{zL(-1)}\Y^{p}(w_{1}, -z)w_{2},
\end{equation}
respectively.

\begin{thm}\label{skew-sym}
The linear maps $\Omega_{+}(\Y)$ and $\Omega_{-}(\Y)$  are twisted intertwining operators
of types ${W_{3}\choose W_{2}\phi_{g_{2}^{-1}}(W_{1})}$ and ${W_{3}\choose \phi_{g_{1}}(W_{2})W_{1}}$,
respectively (recall the definition of $\phi_{g}$ for an automorphism $g$ of $V$ in Section 3).
\end{thm}
\pf
The lower-truncation property and the $L(-1)$-derivative property are easy to verify. We prove 
only the duality property. 

Let $u\in V$, $w_{1}\in W_{1}$, $w_{2}\in W_{2}$
and $w_{3}'\in W_{3}'$.  We first need to give the multivalued analytic functions with preferred branches  in the 
duality property. We shall denote these multivalued analytic functions for $\Omega_{+}(\Y)$ and $\Omega_{-}(\Y)$
by 
$g_{+}(z_{1}, z_{2}; u, w_{1}, w_{2}, w_{3}')$ and 
$g_{-}(z_{1}, z_{2}; u, w_{1}, w_{2}, w_{3}')$, respectively.  Let 
$f(z_1, z_2; u, w_{1}, w_{2}, w_{3}')$ 
be the multivalued analytic function with preferred branch
in the duality property for the twisted intertwining operator $\Y$. Then we can write 
\begin{align*}
f(z_1&-z_{2}, -z_2; u, w_{1}, w_{2}, e^{z_{2}L'(1)}w_{3}')\nn
&=\sum_{i,
j, k, l, m, n= 0}^N a_{ijklmn}(z_1-z_{2})^{r_i}(-z_2)^{s_{j}}z_1^{t_{k}}(\log z_1)^l
(\log z_2)^{m}(\log (z_1 - z_2))^{n},
\end{align*}
where $L'(1)$ is the adjoint operator of $L(-1)$ on $W_{3}$ and is equal to the coefficient of the $x^{-3}$ term in 
$(Y^{g}_{W_{3}})'(\omega, x)$. 
Define
\begin{eqnarray}\label{skew-sym-1}
\lefteqn{g_{\pm}(z_{1}, z_{2}; u, w_{2}, w_{1}, w_{3}')}\nn
&& = \sum_{i,
j, k, l, m, n= 0}^N e^{\pm s_{j}\pi \i}a_{ijklmn}(z_1-z_{2})^{r_i}z_2^{s_{j}}z_1^{t_{k}}(\log (z_1-z_{2}))^l
(\log z_2+\pi \i)^{m}(\log z_1)^{n}.
\end{eqnarray}

When  $|z_{1}|>|z_{2}|>0$, $|z_{1}-z_{2}|>|z_{2}|>0$, $|\arg (z_{1}-z_{2})-\arg z_{1}|<\frac{\pi}{2}$ and 
$\arg z_{2}< \pi$ (for $\Omega_{+}$) or $\arg z_{2}\ge  \pi$ (for $\Omega_{-}$), 
from (\ref{omega3}), the $L(-1)$-derivative property for $Y_{W_{3}}
^{g_{3}}$ and the duality property for $\Y$, 
\begin{eqnarray}\label{skew-sym-0}
\lefteqn{\langle w'_{3}, (Y_{W_{3}}
^{g_{3}})^{p_{1}}(u, z_{1})\Omega_{\pm}(\mathcal{Y})^{p_{2}}(w_{2}, z_{2})w_{1}\rangle}\nn
&&=\langle w'_{3}, (Y_{W_{3}}
^{g_{3}})^{p_{1}}(u, z_{1})e^{z_{2}L(-1)}\mathcal{Y}^{p_{2}}(w_{1}, -z_{2})w_{2}\rangle\nn
&&=\langle e^{z_{2}L'(1)}w'_{3}, (Y_{W_{3}}
^{g_{3}})^{p_{1}}(u, z_{1}-z_{2})\mathcal{Y}^{p_{2}}(w_{1}, -z_{2})w_{2}\rangle
\end{eqnarray}
converges absolutely to 
\begin{eqnarray}\label{skew-sym-2}
\lefteqn{f^{p_{1}, p_{2}, p_{1}}(z_1-z_{2}, -z_2; u, w_{1}, w_{2}, e^{z_{2}L'(1)}w_{3}')}\nn
&&= \sum_{i,
j, k, l, m, n= 0}^Na_{ijklmn}e^{r_{i}l_{p_{1}} (z_1-z_{2})}e^{s_{j}l_{p_{2}}(-z_2)}e^{t_{k}l_{p_{1}}(z_1)} (l_{p_{1}}(z_1-z_{2}))^l
(l_{p_{2}}(-z_2))^{m}(l_{p_{1}}(z_1))^{n}.\nn
\end{eqnarray}
Since $\arg z_{2}< \pi$ (for $\Omega_{+}$) or $\arg z_{2}\ge \pi$ (for $\Omega_{-}$),
$\arg (-z_{2})=\arg z_{2}\pm \pi$ (for $\Omega_{\pm}$). Hence the right-hand side of 
(\ref{skew-sym-2}) is equal to 
\begin{eqnarray}\label{skew-sym-3}
&{\displaystyle \sum_{i,
j, k, l, m, n= 0}^N a_{ijklmn}e^{r_{i}l_{p_{1}} (z_1-z_{2})}
e^{s_{j}(l_{p_{2}}(z_2)\pm \pi \i)}e^{t_{k}l_{p_{1}}(z_1)} (l_{p_{1}}(z_1-z_{2}))^l
(l_{p_{2}}(z_2)\pm \pi \i)^{m}(l_{p_{1}}(z_1))^{n}}\nn
&=g_{\pm}^{p_{1}, p_{2}, p_{1}}(z_{1}, z_{2};  u, w_{2}, w_{1}, w_{3}').&
\end{eqnarray}
But the left-hand side of (\ref{skew-sym-3}) can be expanded in the region $|z_{1}|>|z_{2}|>0$ as a series in 
powers of $e^{l_{p_{1}}(z_{1})}$ and $e^{l_{p_{2}}(z_{2})}$ and in finitely many nonnegative integral powers of 
$l_{p_{1}}(z_{1})$ and $l_{p_{2}}(z_{2})$ such that the real parts of the powers of $e^{l_{p_{1}}(z_{1})}$ is 
bounded from above and the real parts of the powers of $e^{l_{p_{2}}(z_{2})}$ is 
bounded from below. So the left-hand side of (\ref{skew-sym-0})
as a series of the same form that converges to the left-hand side of (\ref{skew-sym-3}) in a smaller region must be 
convergent absolutely  in the larger region given by $|z_{1}|>|z_{2}|>0$ and its sum must be equal 
to (\ref{skew-sym-3}) when $|z_{1}|>|z_{2}|>0$, 
$|\arg (z_{1}-z_{2})-\arg z_{1}|<\frac{\pi}{2}$ and $\arg z_{2}<\pi$ (for $\Omega_{+}$) or
$\arg z_{2}\ge\pi$ (for $\Omega_{-}$). Since both the 
left-hand side of (\ref{skew-sym-0}) and the left-hand side of (\ref{skew-sym-3}) are single-valued 
analytic function in $z_{1}$ and $z_{2}$ with cuts at $z_{1}\in \R_{+}$ and $z_{2}\in \R_{+}$, the fact that
they are equal when $|z_{1}|>|z_{2}|>0$, $|\arg (z_{1}-z_{2})-\arg z_{1}|<\frac{\pi}{2}$ and $\arg z_{2}<\pi$ 
(for $\Omega_{+}$) or $\arg z_{2}\ge\pi$ (for $\Omega_{-}$) means that they are 
equal when $|z_{1}|>|z_{2}|>0$,
$|\arg (z_{1}-z_{2})-\arg z_{1}|<\frac{\pi}{2}$. Thus we have proved that when $|z_{1}|>|z_{2}|>0$ and
$|\arg (z_{1}-z_{2})-\arg z_{1}|<\frac{\pi}{2}$,   the
left-hand side of (\ref{skew-sym-0}) is equal to $g_{\pm}^{p_{1}, p_{2}, p_{1}}(z_{1}, z_{2};  u, w_{2}, w_{1}, w_{3}')$.

When $|z_{2}|>|z_{1}|>0$ and 
$\arg z_{2}\ge \pi$, 
\begin{eqnarray}\label{skew-sym-4}
\lefteqn{\langle w'_{3}, \Omega_{-}(\mathcal{Y})^{p_{2}}(w_{2}, z_{2})(Y_{W_{1}}
^{g_{1}})^{p_{1}}(u, z_{1})w_{1}\rangle}\nn
&&=\langle w'_{3}, e^{z_{2}L(-1)}\mathcal{Y}^{p_{2}}((Y_{W_{1}}
^{g_{1}})^{p_{1}}(u, z_{1})w_{1}, -z_{2})w_{2}\rangle\nn
&&=\langle e^{z_{2}L'(1)}w'_{3}, \mathcal{Y}^{p_{2}}((Y_{W_{1}}
^{g_{1}})^{p_{1}}(u, z_{1})w_{1}, -z_{2})w_{2}\rangle
\end{eqnarray}
converges absolutely and if in addition,  $|\arg (z_{1}-z_{2})-\arg (-z_{2})|<\frac{\pi }{2}$, its sum is equal  to 
\begin{eqnarray}\label{skew-sym-5}
\lefteqn{f^{p_{2}, p_{2}, p_{1}}(z_1-z_{2}, -z_2; u, w_{1}, w_{2}, e^{z_{2}L'(1)}w_{3}')}\nn
&&= \sum_{i,
j, k, l, m, n= 0}^Na_{ijklmn}e^{r_{i}l_{p_{2}} (z_1-z_{2})}e^{s_{j}l_{p_{2}}(-z_2)}e^{t_{k}l_{p_{1}}(z_1)} (l_{p_{2}}(z_1-z_{2}))^l
(l_{p_{2}}(-z_2))^{m}(l_{p_{1}}(z_1))^{n}.\nn
\end{eqnarray}
When $\arg z_{2}\ge \pi$, $\arg (-z_{2})=\arg z_{2}-\pi$. Hence the right-hand side of 
(\ref{skew-sym-5}) is equal to 
\begin{eqnarray}\label{skew-sym-6}
&{\displaystyle \sum_{i,
j, k, l, m, n= 0}^Na_{ijklmn}e^{r_{i}l_{p_{2}} (z_1-z_{2})}e^{s_{j}(l_{p_{2}}(z_2)-\pi \i)}e^{t_{k}l_{p_{1}}(z_1)} (l_{p_{2}}(z_1-z_{2}))^l
(l_{p_{2}}(z_2)-\pi \i)^{m}(l_{p_{1}}(z_1))^{n}}&\nn
&=g_{-}^{p_{1}, p_{2}, p_{2}}(z_{1}, z_{2}; u, w_{1}, w_{2}, w_{3}').&
\end{eqnarray}
The same argument as above shows that 
the left-hand side of (\ref{skew-sym-4}) converges absolutely when $|z_{2}|>|z_{1}|>0$ and 
its sum is equal to (\ref{skew-sym-6}) when $|z_{2}|>|z_{1}|>0$, 
$|\arg (z_{1}-z_{2})-\arg (-z_{2})|<\frac{\pi }{2}$ and $\arg z_{2}\ge \pi$. But when 
$\arg z_{2}\ge \pi$, $\arg (-z_{2})=\arg z_{2}-\pi$. 
Hence in this case, the 
inequality $|\arg (z_{1}-z_{2})-\arg (-z_{2})|<\frac{\pi }{2}$ becomes
$-\frac{3\pi}{2}<\arg (z_{1}-z_{2})-\arg z_{2}<-\frac{\pi}{2}$. 
Also both the left-hand side of 
(\ref{skew-sym-4}) and the left-hand side of (\ref{skew-sym-6})  are single
valued analytic functions in $z_{1}$ and $z_{2}$ with cuts at $z_{1}\in \R_{+}$ and $z_{2}\in \R_{+}$.
Thus the same argument as above shows that
when $|z_{2}|>|z_{1}|>0$ and
$-\frac{3\pi}{2}<\arg (z_{1}-z_{2})-\arg z_{2}<-\frac{\pi}{2}$,  the 
left-hand side of (\ref{skew-sym-4}) is equal to 
$g_{-}^{p_{1}, p_{2}, p_{1}}(z_{1}, z_{2};  u, w_{2}, w_{1}, w_{3}')$.

Next we discuss the iterate of $\Omega_{-}(\Y)$ and the twisted vertex operator
map $\phi_{g_{1}}(Y^{g_{2}}_{W_{2}})$. By Lemma \ref{int-rev-prod-1},
when $|z_{2}|>|z_{1}-z_{2}|>0$ and $\arg z_{2}\ge \pi$, 
\begin{eqnarray}\label{skew-sym-7}
\lefteqn{\langle w'_{3}, \Omega_{-}(\mathcal{Y})^{p_{2}}(\phi_{g_{1}}(Y^{g_{2}}_{W_{2}})
^{p_{12}}(u, z_{1}-z_{2})w_{2}, z_{2})w_{1}\rangle}\nn
&&=\langle w'_{3}, \Omega_{-}(\mathcal{Y})^{p_{2}}((Y_{W_{2}}
^{g_{2}})^{p_{12}}(g_{1}^{-1}u, z_{1}-z_{2})w_{2}, z_{2})w_{1}\rangle\nn
&&=\langle w'_{3}, e^{z_{2}L(-1)}\mathcal{Y}^{p_{2}}(w_{1}, -z_{2})
(Y_{W_{2}}
^{g_{2}})^{p_{12}}(g_{1}^{-1}u, z_{1}-z_{2})w_{2}\rangle\nn
&&=\langle e^{z_{2}L'(1)}w'_{3}, \mathcal{Y}^{p_{2}}(w_{1}, -z_{2})
(Y_{W_{2}}
^{g_{2}})^{p_{12}}(g_{1}^{-1}u, z_{1}-z_{2})w_{2}\rangle
\end{eqnarray}
converges absolutely and if in addition,
$\frac{\pi}{2}<\arg z_{1}-\arg (-z_{2})<\frac{3\pi}{2}$,  its sum
is equal to $f^{p_{12}, p_{2}, p_{2}-1}(z_1-z_{2}, -z_2; g_{1}^{-1}u, w_{1}, w_{2}, e^{z_{2}L'(1)}w_{3}')$.
By (\ref{skew-sym-8.5}),  we have
\begin{eqnarray}\label{skew-sym-9}
\lefteqn{f^{p_{12}, p_{2}, p_{2}-1}(z_1-z_{2}, -z_2; g_{1}^{-1}u, w_{1}, w_{2}, e^{z_{2}L'(1)}w_{3}')}\nn
&&=f^{p_{12}, p_{2}, p_{2}}(z_1-z_{2}, -z_2; u, w_{1}, w_{2}, e^{z_{2}L'(1)}w_{3}')\nn
&&= \sum_{i,
j, k, l, m, n= 0}^Na_{ijklmn}e^{r_{i}l_{p_{12}} (z_1-z_{2})}e^{s_{j}l_{p_{2}}(-z_2)}e^{t_{k}l_{p_{2}}(z_1)} (l_{p_{12}}(z_1-z_{2}))^l
(l_{p_{2}}(-z_2))^{m}(l_{p_{2}}(z_1))^{n}.\nn
\end{eqnarray}
When $\arg z_{1}\ge \pi$, $\arg (-z_{2})=\arg z_{2}-\pi$ and hence the right-hand side of 
(\ref{skew-sym-9}) is equal to 
\begin{eqnarray}\label{skew-sym-10}
&{\displaystyle  \sum_{i,
j, k, l, m, n= 0}^Na_{ijklmn}e^{r_{i}l_{p_{12}} (z_1-z_{2})}e^{s_{j}(l_{p_{2}}(z_2)-\pi \i)}e^{t_{k}l_{p_{2}}(z_1)} (l_{p_{12}}(z_1-z_{2}))^l
(l_{p_{2}}(z_2)-\pi \i)^{m}(l_{p_{2}}(z_1))^{n}}\nn
&=g_{-}^{p_{2}, p_{2}, p_{12}}(z_{1}, z_{2};  u, w_{2}, w_{1}, w_{3}').&
\end{eqnarray}
Thus the left-hand side of (\ref{skew-sym-7}) converges absolutely when $|z_{2}|>|z_{1}-z_{2}|>0$
and its sum is equal to  
$g_{-}^{p_{2}, p_{2}, p_{12}}(z_{1}, z_{2};  u, w_{2}, w_{1}, w_{3}')$ when
$|z_{2}|>|z_{1}-z_{2}|>0$, $\frac{\pi}{2}<\arg z_{1}-\arg (-z_{2})<\frac{3\pi}{2}$ and $\arg z_{2}\ge \pi$.
But when $\arg z_{2}\ge \pi$, $\arg (-z_{2})=\arg z_{2}-\pi$ and thus the inequality 
$\frac{\pi}{2}<\arg z_{1}-\arg (-z_{2})<\frac{3\pi}{2}$ is equivalent to 
$|\arg z_{1}-\arg (z_{2})|<\frac{\pi}{2}$. Then by the same arguments as in the cases above,
we see that when $|z_{2}|>|z_{1}-z_{2}|>0$ and $|\arg z_{1}-\arg z_{2}|<\frac{\pi}{2}$,  the sum of
left-hand side of (\ref{skew-sym-7}) is equal to 
$g_{-}^{p_{2}, p_{2}, p_{12}}(z_{1}, z_{2};  u, w_{2}, w_{1}, w_{3}')$.

We now come back to discuss  $\Omega_{+}(\Y)$.
When $|z_{2}|>|z_{1}|>0$ and 
$\arg z_{2}< \pi$, 
\begin{eqnarray}\label{skew-sym-4-1}
\lefteqn{\langle w'_{3}, \Omega_{+}(\mathcal{Y})^{p_{2}}(w_{2}, z_{2})\phi_{g_{2}^{-1}}(Y_{W_{1}}
^{g_{1}})^{p_{1}}(u, z_{1})w_{1}\rangle}\nn
&&=\langle w'_{3}, e^{z_{2}L(-1)}\mathcal{Y}^{p_{2}}(\phi_{g_{2}^{-1}}(Y_{W_{1}}
^{g_{1}})^{p_{1}}(u, z_{1})w_{1}, -z_{2})w_{2}\rangle\nn
&&=\langle e^{z_{2}L'(1)}w'_{3}, \mathcal{Y}^{p_{2}}((Y_{W_{1}}
^{g_{1}})^{p_{1}}(g_{2}u, z_{1})w_{1}, -z_{2})w_{2}\rangle
\end{eqnarray}
converges absolutely and if in addition, $|\arg (z_{1}-z_{2})-\arg (-z_{2})|<\frac{\pi}{2}$,
its sum is equal to $f^{p_{2}, p_{2}, p_{1}}(z_1-z_{2}, -z_2; g_{2}u, w_{1}, w_{2}, e^{z_{2}L'(1)}w_{3}')$.
By (\ref{skew-sym-8-1.5}),  we have
\begin{eqnarray}\label{skew-sym-9-1}
\lefteqn{f^{p_{2}, p_{2}, p_{1}}(z_1-z_{2}, -z_2; g_{2}u, w_{1}, w_{2}, e^{z_{2}L'(1)}w_{3}')}\nn
&&=f^{p_{2}-1, p_{2}, p_{1}}(z_1-z_{2}, -z_2; u, w_{1}, w_{2}, e^{z_{2}L'(1)}w_{3}')\nn
&&=\sum_{i,
j, k, l, m, n= 0}^Na_{ijklmn}e^{r_{i}l_{p_{2}-1} (z_1-z_{2})}e^{s_{j}l_{p_{2}}(-z_2)}e^{t_{k}l_{p_{1}}(z_1)} (l_{p_{2}-1}(z_1-z_{2}))^l
(l_{p_{2}}(-z_2))^{m}(l_{p_{1}}(z_1))^{n}.\nn
\end{eqnarray}
When $\arg z_{1}< \pi$, $\arg (-z_{2})=\arg z_{2}+\pi$ and hence the right-hand side of 
(\ref{skew-sym-9-1}) is equal to 
\begin{eqnarray}\label{skew-sym-10-1}
&{\displaystyle \sum_{i,
j, k, l, m, n= 0}^Na_{ijklmn}e^{r_{i}l_{p_{2}-1} (z_1-z_{2})}e^{s_{j}(l_{p_{2}}(z_2)+\pi \i)}e^{t_{k}l_{p_{1}}(z_1)} (l_{p_{2}-1}(z_1-z_{2}))^l
(l_{p_{2}}(z_2)+\pi \i)^{m}(l_{p_{1}}(z_1))^{n}}\nn
&=g_{+}^{p_{1}, p_{2}, p_{2}-1}(z_{1}, z_{2};  u, w_{2}, w_{1}, w_{3}').&
\end{eqnarray}
Thus we see that
the left-hand side of (\ref{skew-sym-4-1}) converges absolutely when $|z_{2}|>|z_{1}|>0$ and 
by (\ref{skew-sym-10-1}), 
its sum is equal to $g_{+}^{p_{1}, p_{2}, p_{2}-1}(z_{1}, z_{2};  u, w_{2}, w_{1}, w_{3}')$ when $|z_{2}|>|z_{1}|>0$, 
$|\arg (z_{1}-z_{2})-\arg (-z_{2})|<\frac{\pi}{2}$ and $\arg z_{2}< \pi$. But when 
$\arg z< \pi$, $\arg (-z_{2})=\arg z_{2}+\pi$. 
Hence in this case, the 
inequality $|\arg (z_{1}-z_{2})-\arg (-z_{2})|<\frac{\pi}{2}$ becomes
$\frac{\pi}{2}<\arg (z_{1}-z_{2})-\arg z_{2}<\frac{3\pi}{2}$. 
Also both the left-hand side of 
(\ref{skew-sym-4-1}) and the left-hand side of (\ref{skew-sym-10-1})  are single-valued 
analytic function in $z_{1}$ and $z_{2}$ with cuts at $z_{1}\in \R_{+}$ and $z_{2}\in \R_{+}$.
Thus the same argument as above shows that
when $|z_{2}|>|z_{1}|>0$ and
$\frac{\pi}{2}<\arg (z_{1}-z_{2})-\arg z_{2}<\frac{3\pi}{2}$,  the 
left-hand side of (\ref{skew-sym-4-1}) converges absolutely to 
$g_{+}^{p_{1}, p_{2}, p_{2}-1}(z_{1}, z_{2};  u, w_{2}, w_{1}, w_{3}')$.
Then by Lemma \ref{int-rev-prod-1}, when $|z_{2}|>|z_{1}|>0$ and
$-\frac{3\pi}{2}<\arg (z_{1}-z_{2})-\arg z_{2}<-\frac{\pi}{2}$,  the sum of 
left-hand side of (\ref{skew-sym-4-1}) is equal to 
$g_{+}^{p_{1}, p_{2}, p_{2}}(z_{1}, z_{2};  u, w_{2}, w_{1}, w_{3}')$.

Finally, we discuss the iterate of $\Omega_{+}(\Y)$ and the twisted vertex operator
map $Y^{g_{2}}_{W_{2}}$. When $|z_{2}|>|z_{1}-z_{2}|>0$ and $\arg z_{2}< \pi$, 
\begin{eqnarray}\label{skew-sym-7-1}
\lefteqn{\langle w'_{3}, \Omega_{+}(\mathcal{Y})^{p_{2}}((Y_{W_{2}}
^{g_{2}})^{p_{12}}(u, z_{1}-z_{2})w_{2}, z_{2})w_{1}\rangle}\nn
&&=\langle w'_{3}, e^{z_{2}L(-1)}\mathcal{Y}^{p_{2}}(w_{1}, -z_{2})
(Y_{W_{2}}
^{g_{2}})^{p_{12}}(u, z_{1}-z_{2})w_{2}\rangle\nn
&&=\langle e^{z_{2}L'(1)}w'_{3}, \mathcal{Y}^{p_{2}}(w_{1}, -z_{2})
(Y_{W_{2}}
^{g_{2}})^{p_{12}}(u, z_{1}-z_{2})w_{2}\rangle
\end{eqnarray}
converges absolutely and if in addition, 
$-\frac{3\pi}{2}<\arg z_{1}-\arg (-z_{2})<-\frac{\pi}{2}$, its sum is equal  to 
\begin{eqnarray}\label{skew-sym-5-1}
\lefteqn{f^{p_{12}, p_{2}, p_{2}}(z_1-z_{2}, -z_2; u, w_{1}, w_{2}, e^{z_{2}L'(1)}w_{3}')}\nn
&&= \sum_{i,
j, k, l, m, n= 0}^Na_{ijklmn}e^{r_{i}l_{p_{12}} (z_1-z_{2})}e^{s_{j}l_{p_{2}}(-z_2)}e^{t_{k}l_{p_{2}}(z_1)} (l_{p_{12}}(z_1-z_{2}))^l
(l_{p_{2}}(-z_2))^{m}(l_{p_{2}}(z_1))^{n}.\nn
\end{eqnarray}
When $\arg z_{1}< \pi$, $\arg (-z_{2})=\arg z_{2}+\pi$. Hence the right-hand side of 
(\ref{skew-sym-5-1}) is equal to 
\begin{eqnarray}\label{skew-sym-6-1}
&{\displaystyle   \sum_{i,
j, k, l, m, n= 0}^Na_{ijklmn}e^{r_{i}l_{p_{12}} (z_1-z_{2})}e^{s_{j}(l_{p_{2}}(z_2)+\pi \i)}e^{t_{k}l_{p_{2}}(z_1)} (l_{p_{12}}(z_1-z_{2}))^l
(l_{p_{2}}(z_2)+\pi \i)^{m}(l_{p_{2}}(z_1))^{n}}\nn
&=g_{+}^{p_{1}, p_{2}, p_{2}}(z_{1}, z_{2};  u, w_{2}, w_{1}, w_{3}').&
\end{eqnarray}
The same argument as above shows that 
the left-hand side of (\ref{skew-sym-7-1}) converges absolutely when 
$|z_{2}|>|z_{1}-z_{2}|>0$  and 
its sum is equal to (\ref{skew-sym-6-1}) when $|z_{2}|>|z_{1}-z_{2}|>0$, 
$-\frac{3\pi}{2}<\arg z_{1}-\arg (-z_{2})<-\frac{\pi}{2}$ and $\arg z_{2}< \pi$.  But when 
$\arg z< \pi$, $\arg (-z_{2})=\arg z_{2}+\pi$. 
Hence in this case, the 
inequality $-\frac{3\pi}{2}<\arg z_{1}-\arg (-z_{2})<-\frac{\pi}{2}$  becomes
$|\arg (z_{1}-z_{2})-\arg (-z_{2})|<\frac{\pi }{2}$. 
Also both the left-hand side of 
(\ref{skew-sym-7-1}) and the left-hand side of (\ref{skew-sym-6-1})  are single
valued analytic function in $z_{1}$ and $z_{2}$ with cuts at $z_{1}\in \R_{+}$ and $z_{2}\in \R_{+}$.
Thus the same argument as above shows that
when $|z_{2}|>|z_{1}-z_{2}|>0$, 
$-\frac{3\pi}{2}<\arg z_{1}-\arg (-z_{2})<-\frac{\pi}{2}$,  the sum of the
left-hand side of (\ref{skew-sym-4}) is equal to 
$g_{+}^{p_{1}, p_{2}, p_{1}}(z_{1}, z_{2};  u, w_{2}, w_{1}, w_{3}')$.
\epfv

Let $\mathcal{V}_{W_{1}W_{2}}^{W_{3}}$ be the space of twisted intertwining operators 
of type ${W_{3}\choose W_{1}W_{2}}$. Then we have:

\begin{cor}
The maps $\Omega_{+}: \mathcal{V}_{W_{1}W_{2}}^{W_{3}}\to \mathcal{V}_{W_{2}\phi_{g_{2}^{-1}}(W_{1})}^{W_{3}}$ and 
$\Omega_{-}: \mathcal{V}_{W_{1}W_{2}}^{W_{3}}\to \mathcal{V}_{\phi_{g_{1}}(W_{2})W_{1}}^{W_{3}}$
are linear isomorphisms. In particular, 
$\mathcal{V}_{W_{1}W_{2}}^{W_{3}}$,  $\mathcal{V}_{\phi_{g_{1}}(W_{2})W_{1}}^{W_{3}}$
and $\mathcal{V}_{W_{2}\phi_{g_{2}^{-1}}(W_{1})}^{W_{3}}$
are linearly isomorphic.
\end{cor}
\pf
It is clear that $\Omega_{+}$ and $\Omega_{-}$ are inverse of each other. 
\epfv

The linear isomorphisms $\Omega_{+}$ and $\Omega_{-}$ are called the {\it skew-symmetry isomorphisms}.

\setcounter{equation}{0} \setcounter{thm}{0}

\section{The contragredient isomorphisms}

In this section, we construct what we call the contragredient isomorphisms between the spaces of twisted intertwining
operators of suitable types. These linear isomorphisms  will play an important role in the study of 
rigidity and other related properties of  the still-to-be-constructed $G$-crossed braided
tensor category structure on the category of $g$-twisted $V$-modules for all $g$ in a 
group $G$ of automorphisms of $V$.

Let $g_{1}, g_{2}$ be automorphisms of $V$,  $W_{1}$, $W_{2}$ and $W_{3}$ $g_{1}$-, $g_{2}$- 
and $g_{1}g_{2}$-twisted 
$V$-modules and $\mathcal{Y}$ a twisted intertwining operator 
of type ${W_{3}\choose W_{1}W_{2}}$. 
We define  linear maps
\begin{eqnarray*}
A_{\pm}(\Y): W_{1}\otimes W_{3}'&\to& W_{2}'\{x\}[\log x]\nn
w_{1}\otimes w_{3}'&\mapsto& A_{\pm}(\Y)(w_{1}, x)w_{3}'
\end{eqnarray*}
by
\begin{align}\label{A1}
\langle A_{\pm}(\Y)(w_{1}, x)w_{3}', w_{2}\rangle
&=\langle w_{3}', \Y(e^{xL(1)}e^{\pm \pi \i L(0)}(x^{-L(0)})^{2}w_{1}, x^{-1})w_{2}\rangle
\end{align}
for $w_{1}\in W_{1}$ and $w_{2}\in W_{2}$ and $w_{3}'\in W_{3}'$. 

Let $(W, Y_{W}^{g})$ be a $g$-twisted $V$-module. When $W_{1}=V$, $W_{2}=W_{3}=W$ and 
$\Y=Y_{W}^{g}$, by definition, $A_{+}(Y_{W}^{g})=A_{-}(Y_{W}^{g})=(Y_{W}^{g})'$ (see Section 3). 

Let $L_{W_{1}}^{s}(0)$ be the semisimple part of $L_{W_{1}}(0)$. 
From the definition (\ref{A1}), for $p\in \Z$, $w_{1}\in W_{1}$, $w_{2}\in W_{2}$, $w_{3}'\in W_{3}'$
and $z\in \C^{\times}$, we have 
\begin{align}\label{A2}
\langle A&_{\pm}(\Y)^{p}(w_{1}, z)w_{3}', w_{2}\rangle\nn
&=\langle A_{\pm}(\Y)^{p}(w_{1}, x)w_{3}', w_{2}\rangle\lbar_{x^{n}=e^{nl_{p}(z)},\;\log x=l_{p}(z)}\nn
&=\langle w_{3}', \Y(e^{xL_{W_{1}}(1)}e^{\pm \pi \i L_{W_{1}}(0)}(x^{-L_{W_{1}}(0)})^{2}w_{1}, x^{-1})w_{2}\rangle
\lbar_{x^{n}=e^{nl_{p}(z)},\;\log x=l_{p}(z)}\nn
&=\langle w_{3}', \Y(e^{xL_{W_{1}}(1)}e^{\pm \pi \i L_{W_{1}}(0)}(x^{-L_{W_{1}}^{s}(0)})^{2}\cdot\nn
&\quad\quad\quad\quad\quad\quad\quad\quad\cdot (e^{-(L_{W_{1}}(0)-L_{W_{1}}^{s}(0))\log x})^{2}w_{1}, x^{-1})w_{2}\rangle
\lbar_{x^{n}=e^{nl_{p}(z)},\;\log x=l_{p}(z)}\nn
&=\langle w_{3}', \Y(e^{zL_{W_{1}}(1)}e^{\pm \pi \i L_{W_{1}}(0)}(e^{-l_{p}(z) L_{W_{1}}^{s}(0)})^{2}\cdot\nn
&\quad\quad\quad\quad\quad\quad\quad\quad\cdot(e^{-(L_{W_{1}}(0)-L_{W_{1}}^{s}(0))l_{p}(z)})^{2}w_{1}, y)w_{2}\rangle
\lbar_{y^{n}=e^{-nl_{p}(z)},\;\log y=-l_{p}(z)}\nn
&=\langle w_{3}', \Y(e^{zL_{W_{1}}(1)}e^{\pm \pi \i L_{W_{1}}(0)}e^{-2l_{p}(z) L_{W_{1}}(0)}w_{1}, y)w_{2}\rangle
\lbar_{y^{n}=e^{-nl_{p}(z)},\;\log y=-l_{p}(z)}.
\end{align}
When $\arg z=0$, $\arg z^{-1}=\arg z=0$ and $-l_{p}(z)=l_{-p}(z^{-1})$.
When $\arg z\ne 0$, $\arg z^{-1}=-\arg z+2\pi$ and $-l_{p}(z)=l_{-p-1}(z^{-1})$.
Hence when $\arg z=0$,  the right-hand side of (\ref{A2}) is equal to 
\begin{align}\label{A3}
\langle w_{3}'&, \Y(e^{zL_{W_{1}}(1)}e^{\pm \pi \i L_{W_{1}}(0)}e^{2l_{-p}(z^{-1}) L_{W_{1}}(0)}w_{1}, y)w_{2}\rangle
\lbar_{y^{n}=e^{nl_{-p}(z^{-1})},\;\log y=l_{-p}(z^{-1})}\nn
&=\langle w_{3}', \Y^{-p}(e^{zL_{W_{1}}(1)}e^{\pm \pi \i L_{W_{1}}(0)}e^{2l_{-p}(z^{-1})L_{W_{1}}(0)}
w_{1}, z^{-1})w_{2}\rangle
\end{align}
and  when $\arg z\ne 0$, it is equal to 
\begin{align}\label{A4}
\langle w_{3}'&, \Y(e^{zL_{W_{1}}(1)}e^{\pm \pi \i L_{W_{1}}(0)}e^{2l_{-p-1}(z^{-1}) L_{W_{1}}(0)}
w_{1}, y)w_{2}\rangle
\lbar_{y^{n}=e^{nl_{-p-1}(z^{-1})},\;\log y=l_{-p-1}(z^{-1})}\nn
&=\langle w_{3}', \Y^{-p-1}(e^{zL_{W_{1}}(1)}e^{\pm \pi \i L_{W_{1}}(0)}e^{2l_{-p-1}(z^{-1}) L_{W_{1}}(0)}
w_{1}, z^{-1})w_{2}\rangle.
\end{align}
From (\ref{A2})--(\ref{A4}),  for $w_{1}\in W_{1}$, $w_{2}\in W_{2}$, $w_{3}'\in W_{3}'$
and $z\in \C^{\times}$, 
we have
\begin{equation}\label{A5}
\langle A_{\pm}(\Y)^{p}(w_{1}, z)w_{3}', w_{2}\rangle
=\langle w_{3}', \Y^{-p}(e^{zL_{W_{1}}(1)}e^{\pm \pi \i L_{W_{1}}(0)}e^{2l_{-p}(z^{-1}) L_{W_{1}}(0)}
w_{1}, z^{-1})w_{2}\rangle
\end{equation}
when $\arg z=0$ and 
\begin{equation}\label{A6}
\langle A_{\pm}(\Y)^{p}(w_{1}, z)w_{3}', w_{2}\rangle
=\langle w_{3}', \Y^{-p-1}(e^{zL_{W_{1}}(1)}e^{\pm \pi \i L_{W_{1}}(0)}e^{2l_{-p-1}(z^{-1}) L_{W_{1}}(0)}
w_{1}, z^{-1})w_{2}\rangle
\end{equation}
when $\arg z\ne 0$.

\begin{thm}\label{A}
The linear maps $A_{+}(\Y)$ and $A_{-}(\Y)$  are twisted intertwining operators
of types ${\phi_{g_{1}}(W_{2}')\choose W_{1}W_{3}'}$ and ${W_{2}'\choose W_{1}\phi_{g_{1}^{-1}}(W_{3}')}$, respectively.
\end{thm}
\pf
Just as in the proof of Theorem \ref{skew-sym}, compared with the duality property,  the lower-truncation property 
and the $L_{W_{1}}(-1)$-derivative property can be verified straightforwardly. So here we prove 
only the duality property. 

We first need to give the multivalued analytic functions with preferred branches in the 
duality property. We shall denote these multivalued analytic functions for $A_{+}(\Y)$ and $A_{-}(\Y)$
by 
$h_{+}(z_{1}, z_{2}; u, w_{1}, w_{2}, w_{3}')$ and 
$h_{-}(z_{1}, z_{2}; u, w_{1}, w_{2}, w_{3}')$, respectively.  Let 
$f(z_1, z_2; u, w_{1}, w_{2}, w_{3}')$ 
be the multivalued analytic function with preferred branch
in the duality property for the twisted intertwining operator $\Y$. Then we can write 
\begin{align*}
f(z_1&, z_2; e^{z^{-1}_{1}L_{V}(1)}(-z_{1}^{2})^{L_{V}(0)}u, e^{z^{-1}_{2}L_{W_{1}}(1)}e^{\pm \pi \i L_{W_{1}}(0)}e^{2\log z_{2} L_{W_{1}}(0)}
w_{1}, w_{2}, w_{3}')\nn
&=\sum_{i,
j, k, l, m, n= 0}^N a^{\pm}_{ijklmn}z_1^{r_i}z_2^{s_{j}}(z_1 - z_2)^{t^{k}}(\log z_1)^l
(\log z_2)^{m}(\log (z_1 - z_2))^{n}.
\end{align*}
Define
\begin{align}\label{A-1}
h_{\pm}&(z_{1}, z_{2}; u, w_{2}, w_{1}, w_{3}')\nn
&=  \sum_{i,
j, k, l, m, n= 0}^N e^{\pm t_{k}\pi \i}a^{\pm}_{ijklmn}z_{1}^{-(r_{i}+t_{k})}
z_{2}^{-(s_{j}+t_{k})}(z_{1}-z_{2})^{t_{k}}\cdot\nn
&\quad\quad\quad\quad\quad\quad\quad\quad\quad\quad\cdot  (-\log z_1)^l
(-\log z_2)^{m}(\log (z_1-z_{2})-\log z_{1}-\log z_{2}\pm \pi \i)^{n}.
\end{align}

Let $u\in V$, $w_{1}\in W_{1}$, $w_{2}\in W_{2}$
and $w_{3}'\in W_{3}'$.  We consider $z_{1}, z_{2}\in \C$ satifying
$|z_{2}^{-1}|>|z_{1}^{-1}|>0$ (or equivalently $|z_{1}|>|z_{2}|>0$) and 
 $\arg z_{1}, \arg z_{2}\ne 0$.
Since  $|z_{2}^{-1}|>|z_{1}^{-1}|>0$, 
from (\ref{A6}), $(Y_{W_{2}}
^{g_{2}})'=A_{+}(Y_{W_{2}}
^{g_{2}})$ and the duality property for $\Y$, 
\begin{align}\label{A-0}
\langle \phi&_{g_{1}}((Y_{W_{2}}^{g_{2}})')^{p_{1}}(u, z_{1})A_{+}(\mathcal{Y})^{p_{2}}(w_{1}, z_{2})w_{3}', w_{2}\rangle\nn
&=\langle ((Y_{W_{2}}
^{g_{2}})')^{p_{1}}(g_{1}^{-1}u, z_{1})A_{+}(\mathcal{Y})^{p_{2}}(w_{1}, z_{2})w_{3}', w_{2}\rangle\nn
&=\langle w'_{3}, \mathcal{Y}^{-p_{2}-1}(e^{z_{2}L_{W_{1}}(1)}e^{ \pi \i L_{W_{1}}(0)}e^{2\log (z_{2}^{-1}) L_{W_{1}}(0)}
w_{1}, z_{2}^{-1})\cdot\nn
&\quad\quad\quad\quad\quad\quad\quad\quad\cdot (Y_{W_{2}}
^{g_{2}})^{-p_{1}-1}(e^{z_{1}L_{V}(1)}(-z_{1}^{-2})^{L_{V}(0)}g_{1}^{-1}u, z_{1}^{-1})w_{2}\rangle
\end{align}
converges absolutely and if in addition,  $-\frac{3\pi}{2}<\arg (z_{1}^{-1}-z_{2}^{-1})-\arg z_{2}^{-1}<-\frac{\pi}{2}$, 
its sum is equal  to 
\begin{align}\label{A-2}
f&^{-p_{1}-1, -p_{2}-1, -p_{2}-1}(z_{1}^{-1}, z_{2}^{-1}; \nn
&\quad\quad\quad\quad\quad\quad\quad\quad e^{z_{1}L_{V}(1)}(-z_{1}^{-2})^{L_{V}(0)}g_{1}^{-1}u,
 e^{z_{2}L_{W_{1}}(1)}e^{ \pi \i L_{W_{1}}(0)}e^{2\log (z_{2}^{-1}) L_{W_{1}}(0)}
w_{1}, w_{2}, w_{3}')\nn
&=f^{-p_{1}-1, -p_{2}-1, -p_{2}-1}(z_{1}^{-1}, z_{2}^{-1};\nn
&\quad\quad\quad\quad\quad\quad\quad\quad g_{1}^{-1}e^{z_{1}L_{V}(1)}(-z_{1}^{-2})^{L_{V}(0)}u, 
 e^{z_{2}L_{W_{1}}(1)}e^{ \pi \i L_{W_{1}}(0)}e^{2\log (z_{2}^{-1}) L_{W_{1}}(0)}
w_{1}, w_{2}, w_{3}').
\end{align}
By (\ref{skew-sym-8.5}), (\ref{A-2}) is equal to 
\begin{align}\label{A-2.5}
f&^{-p_{1}-1, -p_{2}-1, -p_{2}}(z_{1}^{-1}, z_{2}^{-1}; e^{z_{1}L_{V}(1)}(-z_{1}^{-2})^{L_{V}(0)}u,
 e^{z_{2}L_{W_{1}}(1)}e^{ \pi \i L_{W_{1}}(0)}e^{2\log (z_{2}^{-1}) L_{W_{1}}(0)}
w_{1}, w_{2}, w_{3}')\nn
&= \sum_{i,
j, k, l, m,n= 0}^N a^{+}_{ijklmn}e^{r_{i}l_{-p_{1}-1} (z_1^{-1})}
e^{s_{j}l_{-p_{2}-1}(z_2^{-1})}e^{t_{k}l_{-p_{2}}(z^{-1}_1-z_{2}^{-1})}\cdot\nn
&\quad\quad\quad\quad\quad\quad\quad\quad\quad\quad\cdot  (l_{-p_{1}-1}(z_1^{-1}))^l
(l_{-p_{2}-1}(z_2^{-1}))^{m}(l_{-p_{2}}(z_1^{-1}-z_{2}^{-1}))^{n}.
\end{align}

Since $\arg z_{1}, \arg z_{2}\ne 0$,  $\arg z_{1}^{-1}=-\arg z_{1}+2\pi$, 
$\arg z_{2}^{-1}=-\arg z_{2}+2\pi$, 
$l_{-p_{1}-1}(z_{1}^{-1})=-l_{p_{1}}(z_{1})$ and 
$l_{-p_{2}-1}(z_{2}^{-1})=-l_{p_{2}}(z_{2})$. Thus 
we also have
$$\arg (z_{1}^{-1}-z_{2}^{-1})=\arg \left(\frac{z_{1}-z_{2}}{z_{1}(-z_{2})}\right).$$
Since $0\le \arg z<2\pi$ for any $z\in \C^{\times}$ and 
$$-\frac{3\pi}{2}<\arg (z_{1}^{-1}-z_{2}^{-1})-\arg z_{2}^{-1}
=\arg \left(\frac{z_{1}-z_{2}}{z_{1}(-z_{2})}\right)-\arg z_{2}^{-1}<-\frac{\pi}{2},$$
we have
\begin{equation}\label{A-3}
\arg \left(\frac{z_{1}-z_{2}}{z_{1}(-z_{2})}\right)=\arg (z_{1}-z_{2})-\arg z_{1}-\arg z_{2}
+(2q+1)\pi
\end{equation}
with $q=-1$ when $-2\pi <\arg (z_{1}-z_{2})-\arg z_{1}<-\frac{3\pi}{2}$, with $q=0$ when 
$|\arg (z_{1}-z_{2})-\arg z_{1}|<\frac{\pi}{2}$ and with $q=1$ when 
$\frac{3\pi}{2} <\arg (z_{1}-z_{2})-\arg z_{1}<2\pi$. From (\ref{A-3}), we obtain 
\begin{align}\label{A-4}
l_{-p_{2}}&(z^{-1}_1-z_{2}^{-1})\nn
&=l_{-p_{2}}\left(\frac{z_{1}-z_{2}}{z_{1}(-z_{2})}\right)\nn
&=\log |z_{1}-z_{2}|-\log |z_{1}|-\log |z_{2}|+\arg (z_{1}-z_{2})\i-\arg z_{1}\i-\arg z_{2}\i\nn
&\quad +(2q+1)\pi \i
+2(-p_{2})\pi \i\nn
&=l_{p_{1}+q}(z_{1}-z_{2})-l_{p_{1}}(z_{1})-l_{p_{2}}(z_{2})+ \pi \i.
\end{align}
From  $l_{-p_{1}-1}(z_{1}^{-1})=-l_{p_{1}}(z_{1})$,
$l_{-p_{2}-1}(z_{2}^{-1})=-l_{p_{2}}(z_{2})$ and (\ref{A-4}),
the right-hand side of (\ref{A-2.5}) is equal to
\begin{align}\label{A-5}
&\sum_{i,
j, k, l, m, n= 0}^N e^{\pm t_{k}\pi \i}a^{\pm}_{ijklmn}e^{-r_{i}l_{p_{1}}(z_{1})}
e^{-s_{j}l_{p_{2}}(z_2)}e^{t_{k}(l_{p_{1}+q}(z_{1}-z_{2})-l_{p_{1}}(z_{1})-l_{p_{2}}(z_{2})+ \pi\i)}\cdot\nn
&\quad\quad\quad\quad\quad\quad\quad\quad\quad\quad\cdot  (-l_{p_{1}}(z_{1}))^l
(-l_{p_{2}}(z_{2}))^{m}(l_{p_{1}+q}(z_{1}-z_{2})-l_{p_{1}}(z_{1})-l_{p_{2}}(z_{2})+\pi \i)^{n}\nn
&\quad =h_{+}^{p_{1}, p_{2}, p_{1}+q}(z_{1}, z_{2}; u, w_{2}, w_{1}, w_{3}').
\end{align}
From (\ref{A-0})--(\ref{A-5}), the left-hand side of (\ref{A-0}) converges absolutely when $|z_{1}|>|z_{2}|>0$ and
$\arg z_{1}, \arg z_{2}\ne 0$ and its sum is equal to the branch 
$h_{+}^{p_{1}, p_{2}, p_{1}+q}(z_{1}, z_{2}; u, w_{2}, w_{1}, w_{3}')$ when $|z_{1}|>|z_{2}|>0$,
$-\frac{3\pi}{2}<\arg (z_{1}^{-1}-z_{2}^{-1})-\arg z_{2}^{-1}<-\frac{\pi}{2}$ and
$\arg z_{1}, \arg z_{2}\ne 0$. In the case that $\arg z_{1}=0$ or $\arg z_{2}= 0$, we can also prove similarly that 
when $|z_{1}|>|z_{2}|>0$, the left-hand side of (\ref{A-0}) converges absolutely and if in addition,
$-\frac{3\pi}{2}<\arg (z_{1}^{-1}-z_{2}^{-1})-\arg z_{2}^{-1}<-\frac{\pi}{2}$, its sum is equal to
$h_{+}^{p_{1}, p_{2}, p_{1}+q}(z_{1}, z_{2}; u, w_{2}, w_{1}, w_{3}')$. The main difference, for example, 
in the case $\arg z_{1}=0$ is that $\arg z_{1}^{-1}=\arg z_{1}=0$ instead of $\arg z_{1}^{-1}=-\arg z_{1}+2\pi$.

When $q=-1$,  since the sum of the left-hand side of (\ref{A-0}) is equal to
the branch $h_{+}^{p_{1}, p_{2}, p_{1}-1}(z_{1}, z_{2}; u, w_{2}, w_{1}, w_{3}')$ 
when $|z_{1}|>|z_{2}|>0$ and  $-2\pi <\arg (z_{1}-z_{2})-\arg z_{1}<-\frac{3\pi}{2}$,  by Lemma \ref{int-rev-prod-1},
the left-hand side of (\ref{A-0}) converges absolutely to
$h_{+}^{p_{1}, p_{2}, p_{1}}(z_{1}, z_{2}; u, w_{2}, w_{1}, w_{3}')$ 
when $|z_{1}|>|z_{2}|>0$ and $0<\arg (z_{1}-z_{2})-\arg z_{1}<\frac{\pi}{2}$. When $q=0$, 
the left-hand side of (\ref{A-0}) converges absolutely to
$h_{+}^{p_{1}, p_{2}, p_{1}}(z_{1}, z_{2}; u, w_{2}, w_{1}, w_{3}')$ 
when $|z_{1}|>|z_{2}|>0$ and $|\arg (z_{1}-z_{2})-\arg z_{1}|<\frac{\pi}{2}$. When $q=1$, 
 since the left-hand side of (\ref{A-0}) converges absolutely to
$h_{+}^{p_{1}, p_{2}, p_{1}+1}(z_{1}, z_{2}; u, w_{2}, w_{1}, w_{3}')$ 
when $|z_{1}|>|z_{2}|>0$ and $\frac{3\pi}{2} <\arg (z_{1}-z_{2})-\arg z_{1}<2\pi$,  by Lemma \ref{int-rev-prod-1},
the sum of the left-hand side of (\ref{A-0}) is equal  to
$h_{+}^{p_{1}, p_{2}, p_{1}}(z_{1}, z_{2}; u, w_{2}, w_{1}, w_{3}')$ 
when $|z_{1}|>|z_{2}|>0$ and $-\frac{\pi}{2}<\arg (z_{1}-z_{2})-\arg z_{1}<0$.
Thus we have proved that when $|z_{1}|>|z_{2}|>0$, $|\arg (z_{1}-z_{2})-\arg z_{1}|<\frac{\pi}{2}$, 
the sum of the left-hand side of (\ref{A-0}) is always equal to
$h_{+}^{p_{1}, p_{2}, p_{1}}(z_{1}, z_{2}; u, w_{2}, w_{1}, w_{3}')$.

Next we consider the product of $A_{+}(\Y)$ and the twisted vertex operator
$(Y_{W_{3}}^{g_{2}})'$. Let $u$, $w_{1}$, $w_{2}$
and $w_{3}'$ be the same as above.  
When  $|z_{1}^{-1}|>|z_{2}^{-1}|>0$ (or equivalently $|z_{2}|>|z_{1}|>0$) and
$\arg z_{1}, \arg z_{2}\ne 0$, 
\begin{align}\label{A-6}
\langle &A_{+}(\mathcal{Y})^{p_{2}}(w_{1}, z_{2})((Y_{W_{2}}^{g_{2}})')^{p_{1}}(u, z_{1})w_{3}', w_{2}\rangle\nn
&=\langle w'_{3}, (Y_{W_{2}}
^{g_{2}})^{-p_{1}-1}(e^{z_{1}L_{V}(1)}(-z_{1}^{-2})^{L_{V}(0)}u, z_{1}^{-1})\cdot\nn
&\quad\quad\quad\quad\quad\quad\quad\quad\cdot 
\mathcal{Y}^{-p_{2}-1}(e^{z_{2}L_{W_{1}}(1)}e^{ \pi \i L_{W_{1}}(0)}e^{2l_{-p_{2}-1}(z_{2}^{-1}) L_{W_{1}}(0)}
w_{1}, z_{2}^{-1}) w_{2}\rangle
\end{align}
converges absolutely and if in addition, $|\arg (z_{1}^{-1}-z_{2}^{-1})-\arg z_{1}^{-1}|<\frac{\pi}{2}$, 
its sum is equal to 
\begin{align}\label{A-7}
f&^{-p_{1}-1, -p_{2}-1, -p_{1}-1}(z_{1}^{-1}, z_{2}^{-1}; \nn
&\quad\quad\quad\quad\quad\quad\quad\quad e^{z_{1}L_{V}(1)}(-z_{1}^{-2})^{L_{V}(0)}u,
 e^{z_{2}L_{W_{1}}(1)}e^{ \pi \i L_{W_{1}}(0)}e^{2l_{-p_{2}-1}(z_{2}^{-1}) L_{W_{1}}(0)}
w_{1}, w_{2}, w_{3}')\nn
&= \sum_{i,
j, k, l, m, n= 0}^N a^{+}_{ijklmn}e^{r_{i}l_{-p_{1}-1} (z_1^{-1})}
e^{s_{j}l_{-p_{2}-1}(z_2^{-1})}e^{t_{k}l_{-p_{1}-1}(z^{-1}_1-z_{2}^{-1})}\cdot\nn
&\quad\quad\quad\quad\quad\quad\quad\quad\quad\quad\cdot  (l_{-p_{1}-1}(z_1^{-1}))^l
(l_{-p_{2}-1}(z_2^{-1}))^{m}(l_{-p_{1}-1}(z_1^{-1}-z_{2}^{-1}))^{n}.
\end{align}

Again we have $\arg z_{1}^{-1}=-\arg z_{1}+2\pi$, 
$\arg z_{2}^{-1}=-\arg z_{2}+2\pi$,
$l_{-p_{1}-1}(z_{1}^{-1})=-l_{p_{1}}(z_{1})$ and 
$l_{-p_{2}-1}(z_{2}^{-1})=-l_{p_{2}}(z_{2})$. 
Since $0\le \arg z<2\pi$ for any $z\in \C^{\times}$ and 
$$-\frac{\pi}{2}<\arg (z_{1}^{-1}-z_{2}^{-1})-\arg z_{1}^{-1}
=\arg \left(\frac{z_{1}-z_{2}}{z_{1}(-z_{2})}\right)-\arg z_{1}^{-1}<\frac{\pi}{2},$$
we have
\begin{equation}\label{A-8}
\arg (z_{1}^{-1}-z_{2}^{-1})=\arg \left(\frac{z_{1}-z_{2}}{z_{1}(-z_{2})}\right)=\arg (z_{1}-z_{2})-\arg z_{1}-\arg z_{2}
+(2q+1)\pi
\end{equation}
with $q=0$ when $\frac{\pi}{2}<\arg (z_{1}-z_{2})-\arg z_{1}<\frac{3\pi}{2}$ and with $q=1$ when 
$-\frac{3\pi}{2}<\arg (z_{1}-z_{2})-\arg z_{1}<-\frac{\pi}{2}$. From (\ref{A-8}) and by the same calculation
as in (\ref{A-4}), we obtain 
\begin{equation}\label{A-9}
l_{-p_{1}-1}(z^{-1}_1-z_{2}^{-1})=l_{p_{2}+q-1}(z_{1}-z_{2})-l_{p_{1}}(z_{1})-l_{p_{2}}(z_{2})+\pi\i.
\end{equation}
 From  $l_{-p_{1}-1}(z_{1}^{-1})=-l_{p_{1}}(z_{1})$,
$l_{-p_{2}-1}(z_{2}^{-1})=-l_{p_{2}}(z_{2})$ and (\ref{A-9}),
the right-hand side of (\ref{A-7}) is equal to
\begin{align}\label{A-10}
\sum_{i,
j, k, l, m, n= 0}^N& a^{+}_{ijklmn}e^{-r_{i}l_{p_{1}}(z_{1})}
e^{-s_{j}l_{p_{2}}(z_{2})}e^{t_{k}(l_{p_{2}+q-1}(z_{1}-z_{2})-l_{p_{1}}(z_{1})-l_{p_{2}}(z_{2})+\pi\i)}\cdot\nn
&\quad\quad\quad\quad\quad\quad\cdot  (-l_{p_{1}}(z_{1}))^l
(-l_{p_{2}}(z_{2}))^{m}(l_{p_{2}+q-1}(z_{1}-z_{2})-l_{p_{1}}(z_{1})-l_{p_{2}}(z_{2})+ \pi\i)^{n}\nn
&=h_{+}^{p_{1}, p_{2}, p_{2}+q-1}(z_{1}, z_{2}; u, w_{2}, w_{1}, w_{3}').
\end{align}
From (\ref{A-6})--(\ref{A-10}), the left-hand side of (\ref{A-6}) converges absolutely when $|z_{2}|>|z_{1}|>0$
and $\arg z_{1}, \arg z_{2}\ne 0$ and if in addition, $|\arg (z_{1}^{-1}-z_{2}^{-1})-\arg z_{1}^{-1}|<\frac{\pi}{2}$, 
its sum is equal to 
$h_{+}^{p_{1}, p_{2}, p_{2}+q-1}(z_{1}, z_{2}; u, w_{2}, w_{1}, w_{3}')$. In the case that $\arg z_{1}=0$ or 
$\arg z_{2}= 0$, we can also prove similarly that 
when $|z_{1}|>|z_{2}|>0$, the left-hand side of (\ref{A-6}) converges absolutely and if in addition,
$-\frac{3\pi}{2}<\arg (z_{1}^{-1}-z_{2}^{-1})-\arg z_{2}^{-1}<-\frac{\pi}{2}$, its sum is equal to
$h_{+}^{p_{1}, p_{2}, p_{1}+q}(z_{1}, z_{2}; u, w_{2}, w_{1}, w_{3}')$.

When $q=0$,  since the sum of the left-hand side of (\ref{A-6}) is equal to the branch
$h_{+}^{p_{1}, p_{2}, p_{2}-1}(z_{1}, z_{2}; u, w_{2}, w_{1}, w_{3}')$ 
when $|z_{2}|>|z_{1}|>0$ and $\frac{\pi}{2}<\arg (z_{1}-z_{2})-\arg z_{1}<\frac{3\pi}{2}$,  
by Lemma \ref{int-rev-prod-1},
the left-hand side of (\ref{A-6}) converges absolutely to
$h_{+}^{p_{1}, p_{2}, p_{2}}(z_{1}, z_{2}; u, w_{2}, w_{1}, w_{3}')$ 
when $|z_{2}|>|z_{1}|>0$ and $-\frac{3\pi}{2}<\arg (z_{1}-z_{2})-\arg z_{1}<-\frac{\pi}{2}$. When $q=1$, 
the left-hand side of (\ref{A-6}) converges absolutely to
$h_{+}^{p_{1}, p_{2}, p_{2}}(z_{1}, z_{2}; u, w_{2}, w_{1}, w_{3}')$ 
when $|z_{2}|>|z_{1}|>0$ and $-\frac{3\pi}{2}<\arg (z_{1}-z_{2})-\arg z_{1}<-\frac{\pi}{2}$.
Thus we have proved that when $|z_{1}|>|z_{2}|>0$, the left-hand side of (\ref{A-6}) converges absolutely 
and if in addition, , 
$-\frac{3\pi}{2}<\arg (z_{1}-z_{2})-\arg z_{1}<-\frac{\pi}{2}$, its sum is equal to
$h_{+}^{p_{1}, p_{2}, p_{2}}(z_{1}, z_{2}; u, w_{2}, w_{1}, w_{3}')$. 

Now we discuss $A_{-}(\mathcal{Y})$. 
When  $|z_{2}^{-1}|>|z_{1}^{-1}|>0$ (or equivalently $|z_{1}|>|z_{2}|>0$) and 
$\arg z_{1}, \arg z_{2}\ne 0$, 
\begin{align}\label{A-10-1}
\langle& ((Y_{W_{2}}^{g_{2}})')^{p_{1}}(u, z_{1})A_{-}(\mathcal{Y})^{p_{2}}(w_{1}, z_{2})w_{3}', w_{2}\rangle\nn
&=\langle w'_{3}, \mathcal{Y}^{-p_{2}-1}(e^{z_{2}L_{W_{1}}(1)}e^{- \pi \i L_{W_{1}}(0)}e^{2l_{-p_{2}-1}(z_{2}^{-1}) L_{W_{1}}(0)}
w_{1}, z_{2}^{-1})\cdot\nn
&\quad\quad\quad\quad\quad\quad\quad\quad\cdot (Y_{W_{2}}
^{g_{2}})^{-p_{1}-1}(e^{z_{1}L_{V}(1)}(-z_{1}^{-2})^{L_{V}(0)}u, z_{1}^{-1})w_{2}\rangle
\end{align}
converges absolutely and if in addition, $-\frac{3\pi}{2}<\arg (z_{1}^{-1}-z_{2}^{-1})-\arg z_{2}^{-1}<-\frac{\pi}{2}$, 
its sum is equal  to 
\begin{align}\label{A-10-2}
f&^{-p_{1}-1, -p_{2}-1, -p_{2}-1}(z_{1}^{-1}, z_{2}^{-1}; \nn
&\quad\quad\quad\quad\quad\quad\quad\quad e^{z_{1}L_{V}(1)}(-z_{1}^{-2})^{L_{V}(0)}u,
 e^{z_{2}L_{W_{1}}(1)}e^{- \pi \i L_{W_{1}}(0)}e^{2l_{-p_{2}-1}(z_{2}^{-1}) L_{W_{1}}(0)}
w_{1}, w_{2}, w_{3}')\nn
&= \sum_{i,
j, k, l, m, n= 0}^N a^{-}_{ijklmn}e^{r_{i}l_{-p_{1}-1} (z_1^{-1})}
e^{s_{j}l_{-p_{2}-1}(z_2^{-1})}e^{t_{k}l_{-p_{2}-1}(z^{-1}_1-z_{2}^{-1})}\cdot\nn
&\quad\quad\quad\quad\quad\quad\quad\quad\quad\quad\cdot  (l_{-p_{1}-1}(z_1^{-1}))^l
(l_{-p_{2}-1}(z_2^{-1}))^{m}(l_{-p_{2}-1}(z_1^{-1}-z_{2}^{-1}))^{n}.
\end{align}

As in the case for $A_{+}(\mathcal{Y})$ above, since $\arg z_{1}, \arg z_{2}\ne 0$, 
$l_{-p_{1}-1}(z_{1}^{-1})=-l_{p_{1}}(z_{1})$ and 
$l_{-p_{2}-1}(z_{2}^{-1})=-l_{p_{2}}(z_{2})$. Also, the same argument as in that case
gives
\begin{equation}\label{A-10-4}
l_{-p_{2}-1}(z^{-1}_1-z_{2}^{-1})=l_{p_{1}+q}(z_{1}-z_{2})-l_{p_{1}}(z_{1})-l_{p_{2}}(z_{2})- \pi\i
\end{equation}
with $q=-1$ when $-2\pi <\arg (z_{1}-z_{2})-\arg z_{1}<-\frac{3\pi}{2}$, with $q=0$ when 
$|\arg (z_{1}-z_{2})-\arg z_{1}|<\frac{\pi}{2}$ and with $q=1$ when 
$\frac{3\pi}{2} <\arg (z_{1}-z_{2})-\arg z_{1}<2\pi$. From  $l_{-p_{1}-1}(z_{1}^{-1})=-l_{p_{1}}(z_{1})$,
$l_{-p_{2}-1}(z_{2}^{-1})=-l_{p_{2}}(z_{2})$ and (\ref{A-10-4}),
the right-hand side of (\ref{A-10-2}) is equal to
\begin{align}\label{A-10-5}
\sum_{i,
j, k, l, m, n= 0}^N& a^{-}_{ijklmn}e^{-r_{i}l_{p_{1}}(z_{1})}
e^{-s_{j}l_{p_{2}}(z_{2})}e^{t_{k}(l_{p_{1}+q}(z_{1}-z_{2})-l_{p_{1}}(z_{1})-l_{p_{2}}(z_{2})- \pi\i)}\cdot\nn
&\quad\quad\quad\quad\quad\quad\cdot  (-l_{p_{1}}(z_{1}))^l
(-l_{p_{2}}(z_{2}))^{m}(l_{p_{1}+q}(z_{1}-z_{2})-l_{p_{1}}(z_{1})-l_{p_{2}}(z_{2})- \pi\i)^{n}\nn
&=h_{-}^{p_{1}, p_{2}, p_{1}+q}(z_{1}, z_{2}; u, w_{2}, w_{1}, w_{3}').
\end{align}
From (\ref{A-10-1})--(\ref{A-10-5}) and the same argument as in the case for $A_{+}(\mathcal{Y})$ above, 
the left-hand side of (\ref{A-10-1}) converges absolutely when 
$|z_{1}|>|z_{2}|>0$ and its sum is equal to the branch
$h_{-}^{p_{1}, p_{2}, p_{1}+q}(z_{1}, z_{2}; u, w_{2}, w_{1}, w_{3}')$ when 
$|z_{1}|>|z_{2}|>0$ and $|\arg (z_{1}-z_{2})-\arg z_{1}|<\frac{\pi}{2}$.

Now we consider the product of $A_{-}(\Y)$ and the twisted vertex operator
$\phi_{g_{1}^{-1}}((Y_{W_{3}}^{g_{2}})')$.
When  $|z_{1}^{-1}|>|z_{2}^{-1}|>0$ (or equivalently $|z_{2}|>|z_{1}|>0$) and
$\arg z_{1}, \arg z_{2}\ne 0$, 
\begin{align}\label{A-10-6}
\langle A_{-}&(\mathcal{Y})^{p_{2}}(w_{1}, z_{2})\phi_{g_{1}^{-1}}((Y_{W_{3}}^{g_{2}})')^{p_{1}}(u, z_{1})w_{3}', w_{2}\rangle\nn
&=\langle A_{-}(\mathcal{Y})^{p_{2}}(w_{1}, z_{2})(Y_{W_{3}}^{g_{2}})'^{p_{1}}(g_{1}u, z_{1})w_{3}', w_{2}\rangle\nn
&=\langle w'_{3}, (Y_{W_{2}}
^{g_{2}})^{-p_{1}-1}(e^{z_{1}L_{V}(1)}(-z_{1}^{-2})^{L_{V}(0)}g_{1}u, z_{1}^{-1})\cdot\nn
&\quad\quad\quad\quad\cdot \mathcal{Y}^{-p_{2}-1}(e^{z_{2}L_{W_{1}}(1)}e^{- \pi \i L_{W_{1}}(0)}e^{2l_{-p_{2}-1}(z_{2}^{-1}) L_{W_{1}}(0)}
w_{1}, z_{2}^{-1}) w_{2}\rangle
\end{align}
converges absolutely and if in addition,  $|\arg (z_{1}^{-1}-z_{2}^{-1})-\arg z_{1}^{-1}|<\frac{\pi}{2}$, 
its sum is equal to 
\begin{align}\label{A-10-7}
f&^{-p_{1}-1, -p_{2}-1, -p_{1}-1}(z_{1}^{-1}, z_{2}^{-1}; \nn
&\quad\quad\quad\quad\quad\quad\quad\quad e^{z_{1}L_{V}(1)}(-z_{1}^{-2})^{L_{V}(0)}g_{1}u, 
e^{z_{2}L_{W_{1}}(1)}e^{- \pi \i L_{W_{1}}(0)}e^{2l_{-p_{2}-1}(z_{2}^{-1}) L_{W_{1}}(0)}
w_{1}, w_{2}, w_{3}')\nn
&=f^{-p_{1}-1, -p_{2}-1, -p_{1}-1}(z_{1}^{-1}, z_{2}^{-1}; \nn
&\quad\quad\quad\quad\quad\quad\quad\quad g_{1}e^{z_{1}L_{V}(1)}(-z_{1}^{-2})^{L_{V}(0)}u, 
e^{z_{2}L_{W_{1}}(1)}e^{- \pi \i L_{W_{1}}(0)}e^{2l_{-p_{2}-1}(z_{2}^{-1}) L_{W_{1}}(0)}
w_{1}, w_{2}, w_{3}').
\end{align}
By (\ref{skew-sym-8.5}), (\ref{A-10-7}) is equal to 
\begin{align}\label{A-10-8}
f&^{-p_{1}-1, -p_{2}-1, -p_{1}-2}(z_{1}^{-1}, z_{2}^{-1}; \nn
&\quad\quad\quad\quad\quad\quad\quad\quad e^{z_{1}L_{V}(1)}(-z_{1}^{-2})^{L_{V}(0)}u, 
 e^{z_{2}L_{W_{1}}(1)}e^{- \pi \i L_{W_{1}}(0)}e^{2l_{-p_{2}-1}(z_{2}^{-1}) L_{W_{1}}(0)}
w_{1}, w_{2}, w_{3}')\nn
&=\sum_{i,
j, k, l, m, n= 0}^N a^{-}_{ijklmn}e^{r_{i}l_{-p_{1}-1} (z_1^{-1})}
e^{s_{j}l_{-p_{2}-1}(z_2^{-1})}e^{t_{k}l_{-p_{1}-2}(z^{-1}_1-z_{2}^{-1})}\cdot\nn
&\quad\quad\quad\quad\quad\quad\quad\quad\quad\quad\cdot  (l_{-p_{1}-1}(z_1^{-1}))^l
(l_{-p_{2}-1}(z_2^{-1}))^{m}(l_{-p_{1}-2}(z_1^{-1}-z_{2}^{-1}))^{n}.
\end{align}

As in the case for $A_{+}(\mathcal{Y})$ above, since $\arg z_{1}, \arg z_{2}\ne 0$, 
$l_{-p_{1}-1}(z_{1}^{-1})=-l_{p_{1}}(z_{1})$ and 
$l_{-p_{2}-1}(z_{2}^{-1})=-l_{p_{2}}(z_{2})$. Also, the same argument as in that case
gives 
\begin{equation}\label{A-10-9}
l_{-p_{1}-2}(z^{-1}_1-z_{2}^{-1})=l_{p_{2}+q-1}(z_{1}-z_{2})-l_{p_{1}}(z_{1})-l_{p_{2}}(z_{2})-\pi\i
\end{equation}
with $q=0$ when $\frac{\pi}{2}<\arg (z_{1}-z_{2})-\arg z_{1}<\frac{3\pi}{2}$ and with $q=1$ when 
$-\frac{3\pi}{2}<\arg (z_{1}-z_{2})-\arg z_{1}<-\frac{\pi}{2}$. From  $l_{-p_{1}-1}(z_{1}^{-1})=-l_{p_{1}}(z_{1})$,
$l_{-p_{2}-1}(z_{2}^{-1})=-l_{p_{2}}(z_{2})$ and (\ref{A-10-9}),
the right-hand side of (\ref{A-10-8}) is equal to
\begin{align}\label{A-10-10}
\sum_{i,
j, k, l, m, n= 0}^N& a^{-}_{ijklmn}e^{-r_{i}l_{p_{1}}(z_{1})}
e^{-s_{j}l_{p_{2}}(z_{2})}e^{t_{k}(l_{p_{2}+q-1}(z_{1}-z_{2})-l_{p_{1}}(z_{1})-l_{p_{2}}(z_{2})- \pi\i)}\cdot\nn
&\quad\quad\quad\quad\quad\quad\cdot  (-l_{p_{1}}(z_{1}))^l
(-l_{p_{2}}(z_{2}))^{m}(l_{p_{2}+q-1}(z_{1}-z_{2})-l_{p_{1}}(z_{1})-l_{p_{2}}(z_{2})- \pi\i)^{n}\nn
&=h_{-}^{p_{1}, p_{2}, p_{2}+q-1}(z_{1}, z_{2}; u, w_{2}, w_{1}, w_{3}').
\end{align}
From (\ref{A-10-6})--(\ref{A-10-10}) and the same argument as in the case for $A_{+}(\mathcal{Y})$ above, 
the left-hand side of (\ref{A-10-6}) converges absolutely when 
$|z_{2}|>|z_{1}|>0$ and its sum is equal to the branch
$h_{-}^{p_{1}, p_{2}, p_{1}+q}(z_{1}, z_{2}; u, w_{2}, w_{1}, w_{3}')$ when 
$|z_{2}|>|z_{1}|>0$ and $-\frac{3\pi}{2}<\arg (z_{1}-z_{2})-\arg z_{1}<-\frac{\pi}{2}$.

Finally we study the iterate of $A_{\pm}(\Y)$ and the twisted vertex operator
$Y_{W_{1}}^{g_{1}}$.
When  
$\arg z_{2}\ne 0$, from (\ref{A6}), we have
\begin{align}\label{A-11}
\langle A_{\pm}&(\mathcal{Y})^{p_{2}}((Y_{W_{1}}^{g_{1}})^{p_{12}}(u, z_{1}-z_{2})w_{1}, z_{2})w_{3}', w_{2}\rangle\nn
&=\langle w'_{3},
\mathcal{Y}^{-p_{2}-1}(e^{z_{2}L_{W_{1}}(1)}e^{\pm \pi \i L_{W_{1}}(0)}e^{2l_{-p_{2}-1}(z_{2}^{-1}) L_{W_{1}}(0)}
 (Y_{W_{1}}^{g_{1}})^{p_{12}}(u, z_{1}-z_{2})w_{1}, z_{2}^{-1}) w_{2}\rangle.
\end{align}
Note that (3.61) and (3.62)  in \cite{HLZ2} with $\Y$ replaced by $Y_{W_{1}}^{g_{1}}$ still holds
since $Y_{W_{1}}^{g_{1}}$
is a (logarithmic) intertwining operator among $V^{\langle g_{1}, g_{2}, g_{1}g_{2}\rangle}$-modules. Using these formulas,
we obtain the formulas for $Y_{W_{1}}^{g_{1}}(u,x)$  ($u\in V$)  conjugated by $e^{x_2L_{W_{1}}(1)}$
and $y^{L_{W_{1}}(0)}$.  Using these formulas,  we obtain
\begin{align}\label{A-11-0}
e^{x_2L_{W_{1}}(1)}&y^{L_{W_{1}}(0)}Y_{W_{1}}^{g_{1}}(u,x_0)\nn
&=e^{x_2L_{W_{1}}(1)}Y_{W_{1}}^{g_{1}}(y^{L_{V}(0)}u, yx_0)y^{L_{W_{1}}(0)}\nn
&=Y_{W_{1}}^{g_{1}}(e^{x_{2}(1-x_{2}yx_{0})L_{V}(1)}(1-x_{2}yx_{0})^{-2L_{V}(0)}y^{L_{V}(0)}u, yx_0(1-x_{2}yx_{0})^{-1})
e^{x_2L_{W_{1}}(1)}y^{L_{W_{1}}(0)}.
\end{align}
Substituting $e^{\pm n\pi \i}x_{2}^{-2n}$ and $\pm \pi \i-2\log x_{2}$ for $y^{n}$ and $\log y$, respectively, 
in (\ref{A-11-0}), we obtain
\begin{align}\label{A-11-0-1}
e^{x_2L_{W_{1}}(1)}&e^{\pm \pi \i L_{W_{1}}(0)}
(x_2^{-L_{W_{1}}(0)})^2 Y_{W_{1}}^{g_{1}}(u,x_0)\nn
&=Y_{W_{1}}^{g_{1}}\left(e^{(x_2+x_0)L_{V}(1)}(-(x_2+x_0)^{-2})^{L_{V}(0)}u,
\frac{x x_0}{(x_2+x_0)x_2}\right) \lbar_{x^{n}=e^{\pm n\pi \i},\; \log x=\pm \pi \i}\cdot\nn
&\quad\quad\cdot e^{x_2L_{W_{1}}(1)}e^{\pm \pi \i L_{W_{1}}(0)}
(x_2^{-L_{W_{1}}(0)})^2
\end{align}
Substituting $e^{nl_{p_{12}}(z_{1}-z_{2})}$, $l_{p_{12}}(z_{1}-z_{2})$, $e^{nl_{p_{2}}(z_{2})}$ 
and $l_{p_{2}}(z_{2})$ for $x^{n}$, $\log x_{0}$, $x_{2}^{n}$ and $\log x_{2}$, respectively, in 
(\ref{A-11-0-1}) and then applying the resulting equality to $w_{1}$, we obtain  in the region $|z_{2}|>|z_{1}-z_{2}|>0$
\begin{align}\label{A-11-1}
e&^{z_{2}L_{W_{1}}(1)}e^{\pm \pi \i L_{W_{1}}(0)}e^{2l_{-p_{2}-1}(z_{2}^{-1}) L_{W_{1}}(0)}
 Y_{W_{1}}^{g_{1}, p_{12}}(u, z_{1}-z_{2})w_{1}\nn
&=Y_{W_{1}}^{g_{1}}\left(e^{z_{1}L_{V}(1)}(-z_{1}^{-2})^{L_{V}(0)}u,
\frac{xx_0}{(x_2+x_0)x_2}\right) \lbar_{\substack{x_{0}^{n}=e^{nl_{p_{12}}(z_{1}-z_{2})},\; 
\log x_{0}=l_{p_{12}}(z_{1}-z_{2}),\; x_{2}^{n}=e^{nl_{p_{2}}(z_{2})}\\ \log x_{2}=l_{p_{2}}(z_{2}),
\;x^{n}=e^{\pm n\pi \i},\;  \log x=\pm \pi \i}}\cdot \nn
&\quad\quad\quad\quad\quad\quad\quad\quad\cdot  e^{z_{2}L_{W_{1}}(1)}e^{\pm \pi \i L_{W_{1}}(0)}e^{2l_{-p_{2}-1}(z_{2}^{-1}) L_{W_{1}}(0)}
w_{1}.
\end{align}

In the region $|z_{2}|>|z_{1}-z_{2}|>0$,  either $\arg (1+\frac{z_{1}-z_{2}}{z_{2}})<\frac{\pi}{2}$
or $\frac{3\pi}{2}<\arg (1+\frac{z_{1}-z_{2}}{z_{2}})$. 
Hence when $|z_{2}|>|z_{1}-z_{2}|>0$, the expansion of $(1+\frac{z_{1}-z_{2}}{z_{2}})^{m}$ for $m\in \C$ as a 
power series in $\frac{z_{1}-z_{2}}{z_{2}}$ is absolutely convergent to 
\begin{equation}\label{A-11-2}
e^{m\log |1+\frac{z_{1}-z_{2}}{z_{2}}|}e^{m(\arg (1+\frac{z_{1}-z_{2}}{z_{2}})\i+2q\pi \i)}
=e^{m\log |z_{1}|}e^{-m\log |z_{2}|}e^{m(\arg (1+\frac{z_{1}-z_{2}}{z_{2}})\i+2q\pi \i)},
\end{equation}
where $q=0$ when $\arg (1+\frac{z_{1}-z_{2}}{z_{2}})<\frac{\pi}{2}$ 
and $q=-1$ when $\frac{3\pi}{2}<\arg (1+\frac{z_{1}-z_{2}}{z_{2}})$.
Also, when $|z_{2}|>|z_{1}-z_{2}|>0$ and $|\arg z_{1}-\arg z_{2}|<\frac{\pi}{2}$, 
\begin{equation}\label{A-11-3}
\arg z_{1}=\arg z_{2}+\arg \left(1+\frac{z_{1}-z_{2}}{z_{2}}\right)+2q\pi,
\end{equation}
where $q=0$ when 
$\arg (1+\frac{z_{1}-z_{2}}{z_{2}})<\frac{\pi}{2}$ and $q=-1$ when 
$\frac{3\pi}{2}<\arg (1+\frac{z_{1}-z_{2}}{z_{2}})$.
By  (\ref{A-11-2}) and (\ref{A-11-3}), we obtain
\begin{align}\label{A-11-4}
x_{2}^{m}&\left(1+\frac{x_0}{x_{2}}\right)^{m}
\lbar_{x_{0}^{n}=e^{nl_{p_{12}}(z_{1}-z_{2})},
\; x_{2}^{n}=e^{nl_{p_{2}}(z_{2})}}\nn
&=e^{ml_{p_{2}}(z_{2})}e^{m\log |z_{1}|}e^{-m\log |z_{2}|}
e^{m((\arg (1+\frac{z_{1}-z_{2}}{z_{2}}))\i+2q\pi \i)}\nn
&=e^{m((\arg z_{2})\i+2p_{2}\pi \i)}e^{m\log |z_{1}|}
e^{m((\arg (1+\frac{z_{1}-z_{2}}{z_{2}})\i)+2q\pi \i)}\nn
&=e^{m((\arg z_{1})\i-(\arg (1+\frac{z_{1}-z_{2}}{z_{2}}))\i-2q\pi \i+2p_{2}\pi \i)}e^{m\log |z_{1}|}
e^{m((\arg (1+\frac{z_{1}-z_{2}}{z_{2}}))\i+2q\pi \i)}\nn
&=e^{ml_{p_{2}}(z_{1})}
\end{align}
for $m\in \C$. Using  (\ref{A-11-4}), we obtain
\begin{align}\label{A-11-5}
&\left(\frac{xx_0}{(x_2+x_0)x_2}\right)^{m}\lbar_{x_{0}^{n}=e^{nl_{p_{12}}(z_{1}-z_{2})},\; x_{2}^{n}
=e^{nl_{p_{2}}(z_{2})},\; x^{n}=e^{\pm n\pi \i}}\nn
&\quad=\left(\frac{e^{\pm m\pi \i}x_0^{m}}
{\displaystyle  \left(x_{2}^{m}\left(1+\frac{x_0}{x_{2}}\right)^{m}\right)x_2^{m}}\right)\lbar_{x_{0}^{n}=e^{nl_{p_{12}}(z_{1}-z_{2})},
\; x_{2}^{n}=e^{nl_{p_{2}}(z_{2})}}\nn
&\quad=e^{\pm m\pi \i}e^{ml_{p_{12}}(z_{1}-z_{2})}e^{-ml_{p_{2}}(z_{1})}e^{-ml_{p_{2}}(z_{2})}\nn
&\quad=e^{m(l_{p_{12}}(z_{1}-z_{2})-l_{p_{2}}(z_{1})-l_{p_{2}}(z_{2})\pm \pi \i)}.
\end{align}
Similarly,  when $|z_{2}|>|z_{1}-z_{2}|>0$, the expansion of $\log (1+\frac{z_{1}-z_{2}}{z_{2}})$ as a 
power series in $\frac{z_{1}-z_{2}}{z_{2}}$ is absolutely convergent to 
\begin{equation}\label{A-11-5.5}
\log \left|1+\frac{z_{1}-z_{2}}{z_{2}}\right|+\arg \left(1+\frac{z_{1}-z_{2}}{z_{2}}\right)\i+2q\pi \i
=\log |z_{1}|-\log |z_{2}|+\arg \left(1+\frac{z_{1}-z_{2}}{z_{2}}\right)\i+2q\pi \i,
\end{equation}
where $q=0$ when $\arg (1+\frac{z_{1}-z_{2}}{z_{2}})<\frac{\pi}{2}$ 
and $q=-1$ when $\frac{3\pi}{2}<\arg (1+\frac{z_{1}-z_{2}}{z_{2}})$.
By  (\ref{A-11-3}) and  (\ref{A-11-5.5}), we obtain
\begin{align}\label{A-11-6}
&\log \left(\frac{xx_0}{(x_2+x_0)x_2}\right)\lbar_{\log x_{0}=l_{p_{12}}(z_{1}-z_{2}),\; \log x_{2}
=l_{p_{2}}(z_{2}),\; \log x=\pm \pi \i}\nn
&\quad=\left(\log x+\log x_{0}-\log x_{2}-\log \left(1+\frac{x_{0}}{x_{2}}\right)-\log x_{2}\right)\lbar_{\log x_{0}=l_{p_{12}}(z_{1}-z_{2}),\; \log x_{2}
=l_{p_{2}}(z_{2}),\; \log x=\pm \pi \i}\nn
&\quad=l_{p_{12}}(z_{1}-z_{2})-l_{p_{2}}(z_{1})-l_{p_{2}}(z_{2})\pm \pi \i.
\end{align}
On the other hand, for any $z_{1}, z_{2}\in \C^{\times}$ such that $z_{1}\ne z_{2}$, there exists
$m\in \Z$ such that 
\begin{align}\label{A-11-7}
l_{p_{12}}&(z_{1}-z_{2})-l_{p_{2}}(z_{1})-l_{p_{2}}(z_{2})\pm \pi \i\nn
&=\log |z_{1}-z_{2}|+(\arg (z_{1}-z_{2}))\i-\log |z_{1}|-(\arg z_{1})\i\nn
&\quad -\log |z_{2}|-(\arg z_{2})\i
+2(p_{12}-2p_{2})\pi \i\pm \pi \i\nn
&=\log \left|\frac{z_{1}-z_{2}}{z_{1}(-z_{2})}\right|+\left(\arg \left(\frac{z_{1}-z_{2}}{z_{1}(-z_{2})}\right)\right)\i
+2(p_{12}-2p_{2})\pi \i+(2m+1)\pi \i\pm \pi \i\nn
&=\log |z_{1}^{-1}-z_{2}^{-1}|+(\arg (z_{1}^{-1}-z_{2}^{-1}))\i
+2\left(p_{12}-2p_{2}+m+\frac{1\pm 1}{2}\right)\pi \i\nn
&=l_{p_{12}-2p_{2}+m+\frac{1\pm 1}{2}}(z_{1}^{-1}-z_{2}^{-1}).
\end{align}

From (\ref{A-11-1}) and (\ref{A-11-5})--(\ref{A-11-7}),  we see that when $|z_{2}|>|z_{1}-z_{2}|>0$ 
and $|\arg z_{1}-\arg z_{2}|<\frac{\pi}{2}$, the right-hand side of (\ref{A-11})
is equal to 
\begin{align}
\langle w'_{3},&
\mathcal{Y}^{-p_{2}-1}((Y_{W_{1}}^{g_{1}})^{p_{12}-2p_{2}+m+\frac{1\pm 1}{2}}(e^{z_{1}L_{V}(1)}(-z_{1}^{-2})^{L_{V}(0)}u, z_{1}^{-1}-z_{2}^{-1})\cdot \nn
&\quad\quad\quad\quad\quad\quad\quad\quad\quad
 \cdot e^{z_{2}L_{W_{1}}(1)}e^{\pm \pi \i L_{W_{1}}(0)}e^{2l_{-p_{2}-1}(z_{2}^{-1}) L_{W_{1}}(0)}
w_{1}, z_{2}^{-1}) w_{2}\rangle,
\end{align}
which, by the duality property for $\Y$, 
converges absolutely to 
\begin{align}\label{A-12}
f&^{-p_{2}-1, -p_{2}-1, p_{12}-2p_{2}+m+\frac{1\pm 1}{2}}(z_{1}^{-1}, z_{2}^{-1}; e^{z_{1}L_{V}(1)}(-z_{1}^{-2})^{L_{V}(0)}u, 
\nn
&\quad\quad\quad\quad\quad\quad\quad\quad\quad\quad\quad\quad
 e^{z_{2}L_{W_{1}}(1)}e^{\pm \pi \i L_{W_{1}}(0)}e^{2l_{-p_{2}-1}(z_{2}^{-1}) L_{W_{1}}(0)}
w_{1}, w_{2}, w_{3}')\nn
&= \sum_{i,
j, k, l, m, n= 0}^N a^{\pm}_{ijklmn}e^{r_{i}l_{-p_{2}-1} (z_1^{-1})}
e^{s_{j}l_{-p_{2}-1}(z_2^{-1})}e^{t_{k}l_{p_{12}-2p_{2}+m+\frac{1\pm 1}{2}}(z^{-1}_1-z_{2}^{-1})}\cdot\nn
&\quad\quad\quad\quad\quad\quad\quad\quad\quad\quad\cdot  (l_{-p_{2}-1}(z_1^{-1}))^l
(l_{-p_{2}-1}(z_2^{-1}))^{m}(l_{p_{12}-2p_{2}+m+\frac{1\pm 1}{2}}(z_1^{-1}-z_{2}^{-1}))^{n}
\end{align}
when $|z_{2}^{-1}|>|z_{1}^{-1}-z_{2}^{-1}|>0$ and  $|\arg z_{1}^{-1}-\arg z_{2}^{-1}|<\frac{\pi}{2}$.
In the case that $\arg z_{1}, \arg z_{2}\ne 0$, $\arg z_{1}^{-1}=-\arg z_{1}+2\pi$, $\arg z_{2}^{-1}=-\arg z_{2}+2\pi$,
$l_{-p_{2}-1}(z_{1}^{-1})=-l_{p_{2}}(z_{1})$ and 
$l_{-p_{2}-1}(z_{2}^{-1})=-l_{p_{2}}(z_{2})$. 
Using these and (\ref{A-11-7}), we see that the right-hand side of (\ref{A-12}) is equal to 
\begin{align}\label{A-13}
&\sum_{i,
j, k, l, m, n= 0}^N a^{\pm}_{ijklmn}e^{-r_{i}l_{p_{2}} (z_1)}
e^{-s_{j}l_{p_{2}}(z_2)}e^{t_{k}(l_{p_{12}}(z_1-z_{2})-l_{p_{2}}(z_{1})-l_{p_{2}}(z_{2})\pm \pi \i)}\cdot\nn
&\quad\quad\quad\quad\quad\quad\quad\quad\quad\quad\cdot  (-l_{p_{2}}(z_1))^l
(-l_{p_{2}}(z_2))^{m}(l_{p_{12}}(z_1-z_{2})-l_{p_{2}}(z_{1})-l_{p_{2}}(z_{2})\pm \pi \i)^{n}\nn
&\quad\quad =h_{\pm}^{p_{2}, p_{2}, p_{12}}(z_{1}, z_{2}; u, w_{2}, w_{1}, w_{3}').
\end{align}

Note that $|z_{2}^{-1}|>|z_{1}^{-1}-z_{2}^{-1}|>0$ is equivalent to $|z_{1}|>|z_{1}-z_{2}|>0$ and 
$|\arg z_{1}^{-1}-\arg z_{2}^{-1}|<\frac{\pi}{2}$ is equivalent to 
$|\arg z_{1}-\arg z_{2}|<\frac{\pi}{2}$. 
Thus we have proved that 
the right-hand side of (\ref{A-11}) is absolutely convergent to 
$h_{\pm}^{p_{2}, p_{2}, p_{12}}(z_{1}, z_{2}; u, w_{2}, w_{1}, w_{3}')$
in the region given by $|z_{1}|, |z_{2}|>|z_{1}-z_{2}|>0$, $|\arg z_{1}-\arg z_{2}|<\frac{\pi}{2}$
and $\arg z_{1}, \arg z_{2}\ne 0$ ($|z_{1}|>|z_{1}-z_{2}|>0$ and $|z_{2}|>|z_{1}-z_{2}|>0$ are both needed in
the proof above). Note that the left-hand side of (\ref{A-11}) is a series 
in (complex) powers of $e^{l_{p_{12}}(z_{1}-z_{2})}$ and $e^{l_{p_{2}}(z_{2})}$
and in nonnegative integral powers of $l_{p_{12}}(z_{1}-z_{2})$
and $l_{p_{2}}(z_{2})$ with 
finitely many negative powers in $e^{l_{p_{12}}(z_{1}-z_{2})}$ and finitely many positive powers of 
$e^{l_{p_{2}}(z_{2})}$.
Since $h_{\pm}^{p_{2}, p_{2}, p_{12}}(z_{1}, z_{2}; u, w_{2}, w_{1}, w_{3}')$
can be expanded uniquely as such a series
in the region  $|z_{2}|>|z_{1}-z_{2}|>0$, 
the left-hand side of (\ref{A-11}) must be absolutely convergent when $|z_{2}|>|z_{1}-z_{2}|>0$
and if in addition, $|\arg z_{1}-\arg z_{1}|<\frac{\pi}{2}$, its sum is equal to 
$h_{\pm}^{p_{2}, p_{2}, p_{12}}(z_{1}, z_{2}; u, w_{2}, w_{1}, w_{3}')$.
\epfv

Just as in the skew-symmetry case, we have the following immediate consequence:

\begin{cor}
The maps $A_{+}: \mathcal{V}_{W_{1}W_{2}}^{W_{3}}\to \mathcal{V}_{W_{1}W_{3}'}^{\phi_{g_{1}}(W_{2}')}$ and 
$A_{-}: \mathcal{V}_{W_{1}W_{2}}^{W_{3}}\to \mathcal{V}_{W_{1}\phi_{g_{1}^{-1}}(W_{3}')}^{W_{2}'}$
are linear isomorphisms. In particular, 
$\mathcal{V}_{W_{1}W_{2}}^{W_{3}}$,  $\mathcal{V}_{W_{1}W_{3}'}^{\phi_{g_{1}}(W_{2}')}$
and $\mathcal{V}_{W_{1}\phi_{g_{1}^{-1}}(W_{3}')}^{W_{2}'}$
are linearly isomorphic.
\end{cor}
\pf
It is clear that $A_{+}$ and $A_{-}$ are inverse of each other. 
\epfv

The linear isomorphisms $A_{+}$ and $A_{-}$ are called the {\it contragredient isomorphisms}.

\noindent {\small \sc Department of Mathematics, Rutgers University,
110 Frelinghuysen Rd., Piscataway, NJ 08854-8019}

\noindent {\em E-mail address}: yzhuang@math.rutgers.edu

\end{document}